\documentclass[final,hidelinks,onefignum,onetabnum]{siamart220329}

\usepackage{lipsum}
\usepackage{amsfonts}
\usepackage{epstopdf}
\usepackage{url}
\usepackage{afterpage}
\ifpdf
  \DeclareGraphicsExtensions{.eps,.pdf,.png,.jpg}
\else
  \DeclareGraphicsExtensions{.eps}
\fi

\usepackage[normalem]{ulem}
\usepackage[labelfont=bf, size=small]{caption}
\usepackage[labelformat=empty, size=small]{subcaption}

\usepackage{microtype, fancyhdr,
  siunitx, graphicx, hyperref, epstopdf, tikz, adjustbox, mathtools, float,
  anyfontsize, relsize}
\usepackage{bm, amssymb, amsmath, mathrsfs}
\usepackage{etoolbox}
\usepackage[algo2e, ruled, linesnumbered, resetcount]{algorithm2e}
\SetFuncSty{textsc}
\SetKwProg{Function}{}{:}{}
\SetKwComment{Comment}{$\triangleright$\ }{}

\usepackage{enumitem}
\setlist[itemize]{leftmargin=2em, topsep=0.5\baselineskip, itemsep=0.5\baselineskip}

\newcommand{\Z}{\mathbb{Z}}    % Integers
\newcommand{\R}{\mathbb{R}}    % Real numbers
\newcommand{\C}{\mathbb{C}}    % Complex numbers
\newcommand{\bO}{\mathcal{O}}  % "big O"
\newcommand{\pluseq}{\ \mathsmaller{\mathrel{+}=}\ } % plus equals
\newcommand{\abs}[1]{\left|#1\right|}      % absolute value
\newcommand{\ind}[1]{\bm{1}_{\set{#1}}}    % Indicator function
\newcommand{\dif}[1]{\mathop{}\!\mathrm{d}{#1}}   % d for dx in integral
\newcommand{\norm}[2][]{\left\Vert#2\right\Vert_{#1}} % norm
\newcommand{\set}[1]{\left\{ #1\right\}}              % set notation
\newcommand{\flr}[1]{{\left\lfloor #1 \right\rfloor}} % floor function
\renewcommand{\epsilon}{\varepsilon}

\DeclareMathOperator*{\argmax}{arg\,max}   % argmax
\let\Re\relax \DeclareMathOperator*{\Re}{Re} 
\let\Im\relax \DeclareMathOperator*{\Im}{Im}
\numberwithin{equation}{section}

\newcommand{\mtx}[1]{\bm{\mathsf{#1}}}

\newcommand{\vct}[1]{\bm{\mathsf{#1}}}

\newcommand{\Bloc}{B^{\textsc{\tiny \hspace*{-1pt} loc}}_{\nu,L}} 
\newcommand{\Basy}{B^{\textsc{\tiny \hspace*{-1pt} asy}}_{\nu,M}} 

\def\spacingset#1{\renewcommand{\baselinestretch}%
{#1}\small\normalsize} \spacingset{1}

\newsiamremark{remark}{Remark}
\newsiamremark{hypothesis}{Hypothesis}
\crefname{hypothesis}{Hypothesis}{Hypotheses}
\newsiamthm{claim}{Claim}

\headers{A Nonuniform Fast Hankel Transform}{Paul G. Beckman and Michael O'Neil}

\title{A Nonuniform Fast Hankel Transform\thanks{Submitted to the editors DATE.
    \funding{P. G. Beckman was partially supported by the Office of Naval
      Research under award \#N00014-21-1-2383 and by the U.S.  Department of
      Energy, Office of Science, Office of Advanced Scientific Computing
      Research, Department of Energy Computational Science Graduate Fellowship
      under Award Number DE-SC0022158. M. O'Neil was partially supported by the
      Office of Naval Research under award \#N00014-21-1-2383.}  } }

\author{Paul G. Beckman\thanks{Courant Institute, New York University, New York, NY\\
  (\email{paul.beckman@cims.nyu.edu}, \url{https://cims.nyu.edu/\~pgb8409}).}
\and Michael O'Neil\thanks{Courant Institute, New York University, New York, NY\\
  (\email{oneil@cims.nyu.edu}, \url{https://cims.nyu.edu/\~oneil}).}
}

\usepackage{amsopn}
\DeclareMathOperator{\diag}{diag}

\ifpdf
\hypersetup{
  pdftitle={A Nonuniform Fast  Hankel Transform}
  pdfauthor={P. G. Beckman and M. O'Neil}
}
\fi

\begin{document}

\maketitle

% REQUIRED
\begin{abstract}
  We describe a fast algorithm for computing discrete Hankel transforms of
moderate orders from $n$ nonuniform points to $m$ nonuniform frequencies in
$\bO\big((m+n)\log\min(n,m)\big)$ operations. Our approach combines local and
asymptotic Bessel function expansions with nonuniform fast Fourier transforms.
The order of each expansion is adjusted automatically according to error
analysis to obtain any desired precision $\epsilon$. Several numerical examples
are provided which demonstrate the speed and accuracy of the algorithm in
multiple regimes and applications.

%%% Local Variables: %% mode: latex %% TeX-master: "../main" %% End:

\end{abstract}

% REQUIRED

% REQUIRED

\section{Introduction} 

The fast Fourier transform (FFT) has revolutionized a wide range of applications
across mathematics, statistics, and the physical sciences by enabling signal
processing and Fourier analysis tasks to be performed using a computational cost
which scales quasi-linearly with the number of data points~$n$. However, the FFT
requires that the input signal be sampled at equispaced points in time and that
the desired output frequencies are equispaced on the integers. These assumptions
are frequently not met in applications such as adaptive numerical partial
differential equation (PDE)
solvers~\cite{alpert2002adaptive,jiang2023dual,askham2017adaptive,nochetto2009theory},
magnetic resonance
imaging~\cite{greengard2007fast,bondesson2019nonuniform,bronstein2002reconstruction},
and various signal processing
tasks~\cite{alexander2012adaptive,thakur2011synchrosqueezing}. To overcome this
setback, nonuniform FFT (NUFFT) algorithms have been
developed~\cite{dutt1993fast,greengard2004accelerating} which achieve near-FFT
speeds in one dimension, assuming that the distribution of time samples and
frequency outputs is not pathological. In higher dimensions, NUFFTs are less
competitive with standard FFTs, but the computational task at hand is also
significantly harder.

The FFT and NUFFT grew out of a need to perform Fourier transforms in Cartesian
coordinates. However, depending on the particular problem, the relevant
continuous Fourier analysis might be better suited to other coordinate systems.
One such commonly encountered situation is computing the Fourier transform of
radially symmetric functions in dimensions~$d \geq 2$. For example, in two
dimensions the Fourier transform of a function~$f$ is given by
\begin{equation}
  g(\omega_1, \omega_2) = \frac{1}{4\pi^2} \iint_{\R^2} f(x_1, x_2) \, 
  e^{-i(\omega_1 x_1 + \omega_2 x_2)}  \, dx_1 \, dx_2.
\end{equation}
Transforming to polar coordinates~$(\omega_1,\omega_2) \mapsto (\omega,\alpha)$
and~$(x_1,x_2) \mapsto (r,\theta)$ the above expression becomes
\begin{equation}
  \label{eq:ftpolar}
  \begin{aligned}
    g(\omega, \alpha) &= \frac{1}{4\pi^2} \int_0^{2\pi} \int_0^\infty
    f(r,\theta) \, 
    e^{-i \omega r (\cos\alpha \cos\theta + \sin\alpha \sin\theta) } 
    \, r \, dr \, d\theta \\
  &= \frac{1}{4\pi^2} \int_0^{2\pi} \int_0^\infty f(r,\theta) \, e^{-i \omega r \cos(\alpha-\theta) } \, r \, dr \, d\theta.
  \end{aligned}
\end{equation}
Furthermore, if~$f$ is radially symmetric, i.e.~$f(r, \theta) = f(r)$, then the
above transform can be written as
\begin{equation}
  \label{eq:HT}
  \begin{aligned}
  g(\omega,\alpha) &= \frac{1}{4\pi^2} \int_0^\infty f(r) \, r \int_0^{2\pi} 
  e^{-i \omega r \cos(\alpha - \theta) }  \, d\theta \, dr \\
  &= \frac{1}{2\pi} \int_0^\infty f(r) \, J_0(\omega r) \, r \, dr,
  \end{aligned}
\end{equation}
where we have used the integral representation of the zeroth-order Bessel
function~\cite{olver2010nist}
\begin{equation}
  J_0(x) 
  = \frac{1}{\pi} \int_0^\pi \cos \left( x \cos \theta \right) \, d\theta.
\end{equation}
The final integral involving~$J_0$ in equation~\eqref{eq:HT} is known as a
\emph{Hankel Transform} of order 0 --- usually referred to simply as a Hankel
Transform. 

In higher ambient dimensions, the Fourier transform of radially symmetric
functions reduces to a Hankel transform of higher order. Similarly, if the
function~$f$ in~\eqref{eq:ftpolar} has a particular periodic dependence
in~$\theta$ so that $f(r,\theta) = f(r)e^{i\nu\theta}$ with $\nu \in \Z$, then
we have
\begin{equation} \label{eq:FB-integral}
  \begin{aligned}
  g(\omega,\alpha) &= \frac{1}{4\pi^2} \int_0^\infty f(r) \, r \int_0^{2\pi} 
  e^{-i \omega r \cos(\alpha - \theta) } \, e^{i\nu\theta}  \, d\theta \, dr \\
  &= \frac{i^\nu}{2\pi} \int_0^\infty f(r) \, r \, J_\nu(\omega r)  \, dr,
  \end{aligned}
\end{equation}
where, again, we have invoked an integral representation
for~$J_\nu$~\cite{olver2010nist}. 

In order to numerically compute~$g$ in~\eqref{eq:HT} or~\eqref{eq:FB-integral}
at a collection of~$m$ ``frequencies''~$\omega_j$, the Hankel transform must be
discretized using an appropriate quadrature rule with nodes $r_k$ and weights
$w_k$ which depend on the particular class of~$f$ for which the integral is
desired. In general this results in the need for computing
\begin{equation} \label{eq:DHT}
  \begin{aligned}
  g(\omega_j) \approx 
  g_j &:= \sum_{k=1}^n w_k \, f(r_k) \, r_k \, J_\nu(\omega_j r_k) \\
  &\ = \sum_{k=1}^n c_k \, J_\nu(\omega_j r_k)
   \qquad \text{for } j = 1, \ldots, m.
  \end{aligned}
\end{equation}
The above sum will be referred to as the Discrete Hankel Transform (DHT) of
order $\nu$. 

In our motivating example --- computing the continuous Fourier transform --- the
DHT arises from the discretization of the radially symmetric Fourier integral.
The DHT also appears in a wide range of applications including
imaging~\cite{higgins1988hankel, zhao2013fourier, marshall2023fast},
statistics~\cite{lord1954a, genton2002nonparametric}, and separation of
variables methods in partial differential
equations~\cite{bisseling1985fast,ali1999generalized, zhou2022spectral}. In many
such applications, a fully nonuniform DHT is desired, as the relevant
frequencies $\omega_j$ may not be equispaced, and the most efficient quadrature
rule for discretizing (\ref{eq:HT}) may have nodes $r_k$ which are also not
equispaced. 

The algorithm of this work allows for arbitrary selection of the
frequencies~$\omega_j$ and nodes~$r_k$, in contrast to other algorithms which
require some structure to their location (e.g. equispaced or exponentially
distributed). There are a few types of commonly encountered DHTs, all of which
our algorithm can address. Schl\"omilch
expansions~\cite{linton2006schlomilch,townsend2015fast} take
frequencies~$\omega_j = j\pi$. Fourier-Bessel expansions --- often used in
separation of variables calculations for PDEs --- take frequencies~$\omega_j =
\beta_{\nu,j}$, where $\beta_{\nu,j}$ denotes the $j^{th}$ root of $J_\nu$. In
the most restrictive cases~\cite{johnson1987improved}, one fixes both~$\omega_j
= \beta_{\nu,j}$ and~$r_k = \beta_{\nu,k}/\beta_{\nu,k+1}$.

\subsection*{Existing methods}
\label{sec:existing}

A number of methods exist in the literature to evaluate (\ref{eq:HT}) and
(\ref{eq:DHT}). These include series expansion methods
\cite{lord1954b,brunol1977fourier,cavanagh1979numerical}, convolutional
approaches \cite{siegman1977quasi, johansen1979fast, mook1983algorithm,
liu1999nonuniform}, and projection-slice or Abel transform-based methods
\cite{oppenheim1980computation, hansen1985fast, kapur1995algorithm}.
See~\cite{cree1993algorithms} for a review of many of these early computational
approaches. Unfortunately, these existing methods are either not applicable to
the discrete case, require a particular choice of $\omega_j$ or $r_k$ due to the
constraints of interpolation or quadrature subroutines, or suffer from low
accuracy as a result of intermediate approximations. Therefore, extending these
schemes to compute the fully nonuniform DHT with controllable accuracy is not
straightforward.

A notable contribution is~\cite{liu1999nonuniform}, which describes a fully
nonuniform fast Hankel transform. This work takes the popular convolutional
approach, using a change of variables to reformulate the Hankel transform as a
convolution with a known kernel which can be evaluated using the NUFFT. However,
its accuracy is limited by the need for a quadrature rule on the nonuniform
points $r_k$. The authors use an irregular trapezoidal rule for this purpose,
which is not high-order accurate. This method also requires the computation of
the inverse NUFFT using conjugate gradients. For even moderately clustered
points or frequencies, this inverse problem is extremely ill-conditioned, and
thus the number of required iterations can be prohibitive. This method is
therefore suitable for ``quasi-equispaced'' points and frequencies, but is not
tractable in general.

More recently, butterfly algorithms \cite{oneil2010algorithm, li2015butterfly,
pang2020interpolative} were introduced as a broadly applicable methodology for
rapidly computing oscillatory transforms including the nonuniform DHT. However,
these algorithms require a precomputation or factorization stage for each new
set of $\omega_j$ and $r_k$. Such precomputations can, unfortunately, be a
bottleneck for applications in which these evaluation points change with each
iteration or application of the transform. In order to provide a
precomputation-free fast DHT,~\cite{townsend2015fast} employs a combination of
asymptotic expansions and Bessel function identities evaluated using the
equispaced FFT. The resulting scheme is applicable to equispaced or perturbed
``quasi-equispaced'' grids in space and frequency, for example $\omega_j =
\beta_{0,j}$ and $r_k = \beta_{0,k} / \beta_{0,n+1}$.

\subsection*{Novelty of this work}
\label{sec:novelty}

% At a high level, our algorithm can be viewed as a generalization of the one
% described in~\cite{townsend2015fast}. In~\cite{townsend2015fast}, asymptotic
% expansions were used to replace~$J_0$ for various arguments. These asymptotic
% expansions involved trigonometic functions, resulting in a fast algorithm for
% computing the DHT using fast cosine transforms (FCTs) and fast sine transforms
% (FSTs). In order to invoke these fast algorithms, various assumptions
% on~$\omega_j$ and~$r_k$ had to be made.

We describe here a precomputation-free nonuniform fast Hankel transform (NUFHT)
which generalizes~\cite{townsend2015fast} to the fully nonuniform setting in a
number of ways. First, we employ an adaptive partitioning scheme which, for any
choice of $\omega_j$ and $r_k$, subdivides the matrix with
entries~$J_\nu(\omega_j r_k)$ into blocks for which matrix-vector products can
be evaluated efficiently. Second, we use the NUFFT to evaluate asymptotic
expansions for nonuniform $r_k$ and $\omega_j$. Finally, we utilize the low-rank
expansion of $J_\nu$ given in~\cite{wimp1962polynomial} in the local regime
where asymptotic expansions are not applicable. We derive error bounds for this
low-rank expansion, allowing us to choose all approximation parameters
automatically by analysis which guarantees that the resulting error is bounded
by the user-specified tolerance $\epsilon$.

\subsection*{Outline of the paper}

The paper is organized as follows. In Section~\ref{sec:overview} we give a high
level view of our algorithm, omitting technical details. Then in
Section~\ref{sec:approx} we study the local and asymptotic expansions of Bessel
functions which serve as the key building blocks of the algorithm. Afterward, in
Section~\ref{sec:methods}, we provide a detailed description of the algorithm
and its associated complexity. Various numerical examples are provided in
Section~\ref{sec:results}, and we conclude with some additional discussion in
Section~\ref{sec:discussion}.

%%% Local Variables: %% mode: latex %% TeX-master: "../main" %% End:

\section{Overview of the algorithm} \label{sec:overview}

To more concisely describe our approach, we write the DHT (\ref{eq:DHT}) as the
equivalent matrix-vector product with $\mtx{A} \in \R^{m \times n}$
\begin{equation}
    \vct{g} = \mtx{A}\vct{f}, \qquad \mtx{A}(j,k) = J_\nu(\omega_j r_k).
\end{equation}
The matrix $\mtx{A}$ is in general full rank and possesses complex oscillatory
structure. As a result, no straightforward fast algorithm exists to apply the
full matrix $\mtx{A}$ to a vector. However, we design an NUFHT by noting that
certain \textit{blocks} $\mtx{A}(j_0:j_1, k_0:k_1)$ are able to be applied to a
vector rapidly using analytical expansions of the underlying Bessel function
$J_\nu$. 

\begin{figure}[!t]
  \centering
  \begin{subfigure}[b]{0.45\textwidth}
    \includegraphics[width=\textwidth]{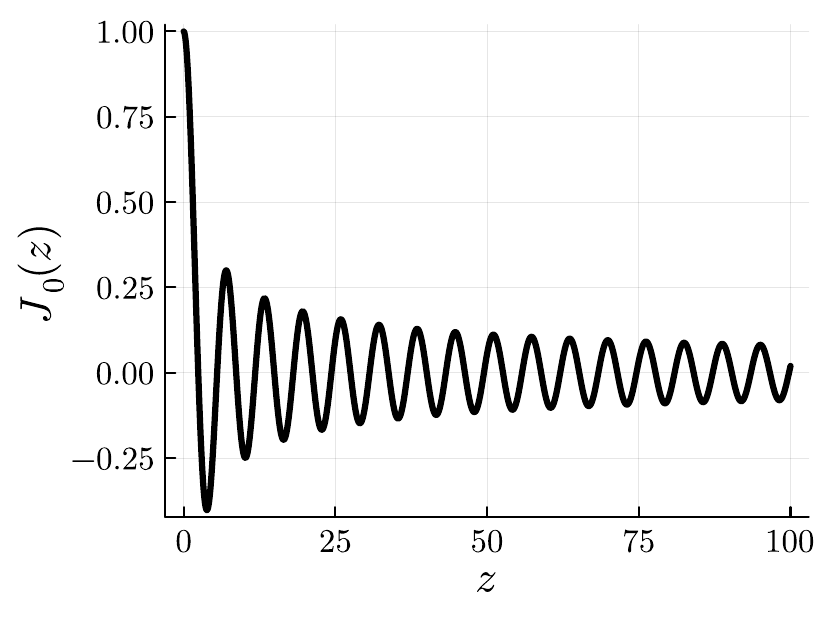}
  \end{subfigure}
  \begin{subfigure}[b]{0.45\textwidth}
    \includegraphics[width=\textwidth]{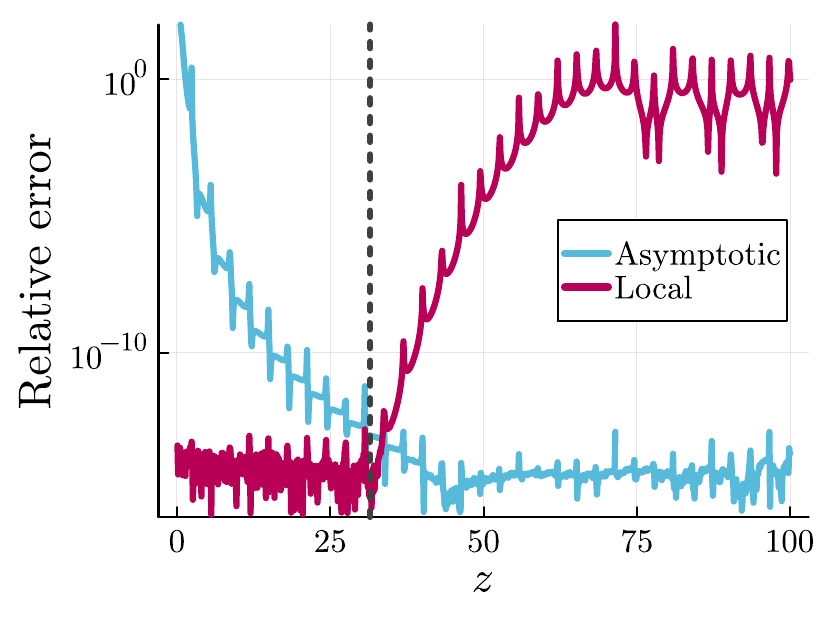}
  \end{subfigure}
  \caption{Bessel function $J_0(z)$ and pointwise relative error in
  approximating $J_0(z)$ using 31-term local and 4-term asymptotic expansions.
  Dotted vertical line shows crossover point where both expansions are accurate
  to $\epsilon = 10^{-12}$.}
  \label{fig:two-expansions}
\end{figure}

When the argument $\omega_j r_k$ is small, $J_\nu$ is smooth and essentially
non-oscillatory, and we use a closed-form local expansion which approximates
$J_\nu$ in terms of Chebyshev polynomials, yielding a low-rank approximation to
various matrix blocks that can be applied to a vector in linear time. When the
argument $\omega_j r_k$ is large, we use a classical asymptotic expansion which
expresses $J_\nu$ as a sum of a small number of decaying sinusoids, and can
therefore be applied to a vector in quasilinear time using the NUFFT.
Figure~\ref{fig:two-expansions} shows the oscillatory behavior of $J_0$, as well
as the accuracy of these local and asymptotic expansions.

\begin{figure}
  \centering
  \newcommand\twa{0.29cm} \newcommand\tw{0.43cm}
  \begin{subfigure}[b]{0.28\textwidth}
    \begin{tikzpicture}
        \draw (0, 0) node[inner sep=0] {\includegraphics[width=0.75\textwidth,
        trim={\twa, \twa, \twa, \twa}, clip]{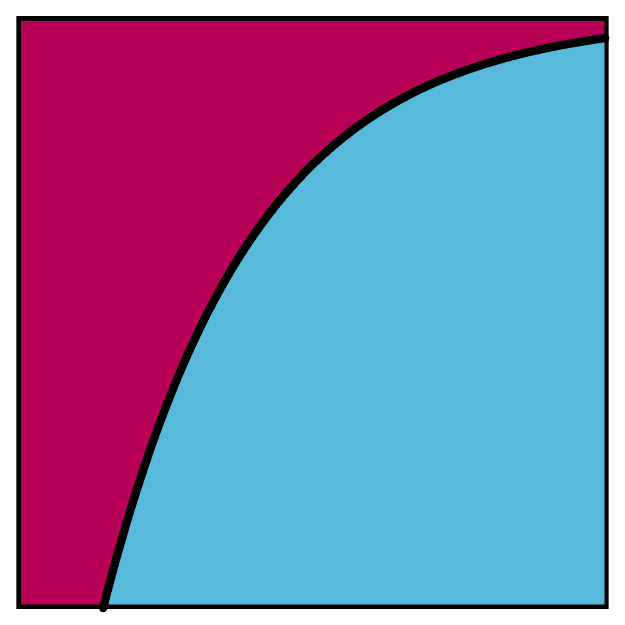}}; \draw
        (-0.8, 0.8) node {\small \textbf{Local}}; \draw (0.3, -0.7) node {\small
        \textbf{Asymptotic}}; \draw (0, 1.6) node {\small $r_1 \ \; < \ \dots \
        < \ \; r_n$}; \draw (-1.7, 0) node[align=center] {\small $\omega_1$
        \\[3pt] $\wedge$ \\[3pt] $\vdots$ \\[3pt] $\wedge$ \\[3pt] $\omega_m$};
    \end{tikzpicture}
    \caption{}
  \end{subfigure}
  \hspace*{0.0\textwidth}
  \begin{subfigure}[b]{0.21\textwidth}
    \includegraphics[width=\textwidth, trim={\tw, \tw, \tw, \tw}, clip]{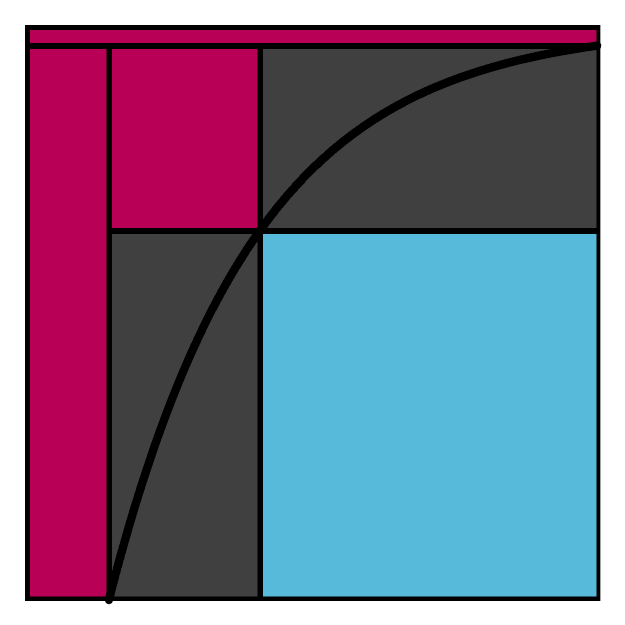}
    \caption{Level 1}
  \end{subfigure}
  \hfill
  \begin{subfigure}[b]{0.21\textwidth}
    \includegraphics[width=\textwidth, trim={\tw, \tw, \tw, \tw}, clip]{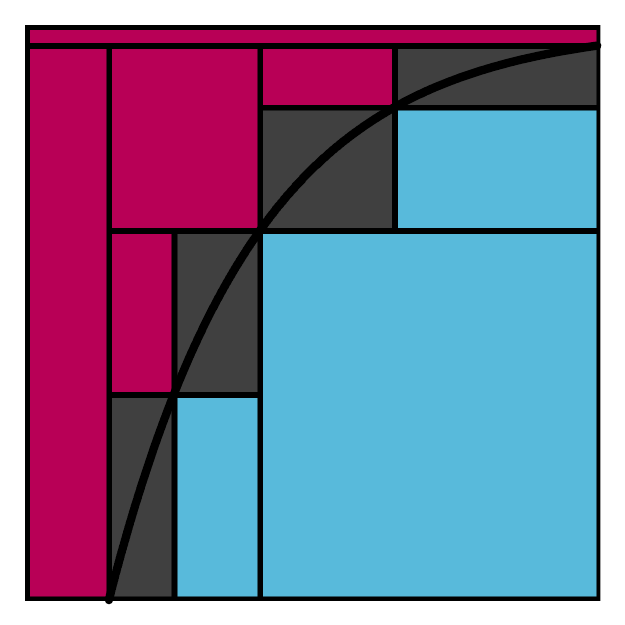}
    \caption{Level 2}
  \end{subfigure}
  \hfill
  \begin{subfigure}[b]{0.21\textwidth}
    \includegraphics[width=\textwidth, trim={\tw, \tw, \tw, \tw}, clip]{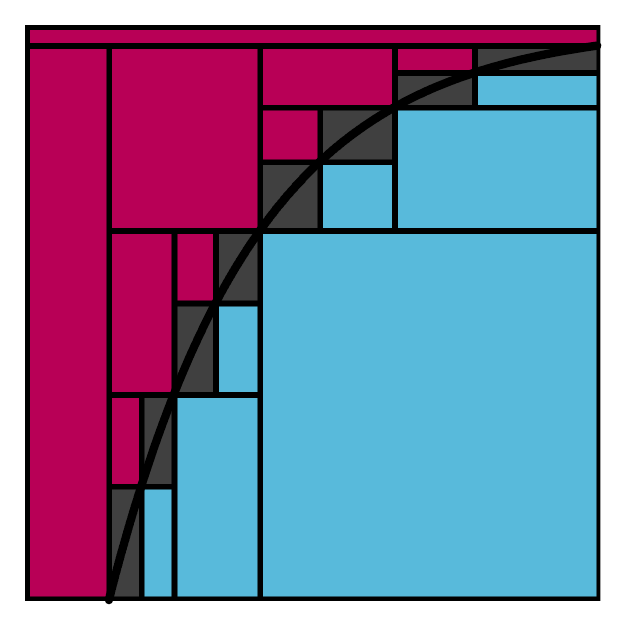}
    \caption{Level 3}
  \end{subfigure}
  \caption{Splitting of Hankel transform matrix $\mtx{A}$ along the curve
  $\omega r = z$ into local and asymptotic regions. Adaptive subdivision of
  $\mtx{A}$ into corresponding local (red), asymptotic (blue), and mixed (gray)
  sub-blocks at various levels.}
  \label{fig:subdivide}
\end{figure}
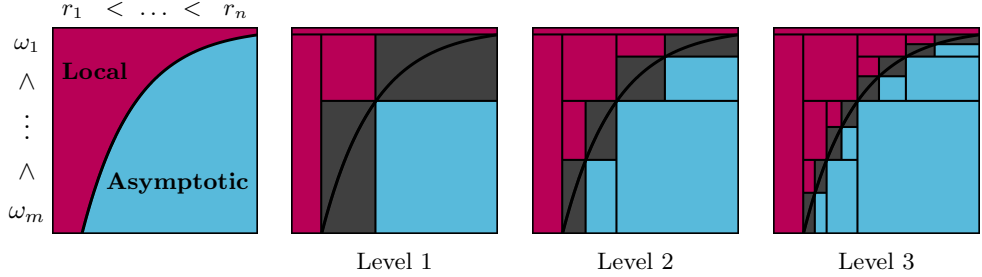

By analyzing the error in these two expansions, we can choose a crossover point
$z$ such that an $L$-term local expansion and an $M$-term asymptotic expansion
are both guaranteed to be accurate to the desired tolerance $\epsilon$ in the
regions $\omega_j r_k \leq z$ and $\omega_j r_k > z$ respectively. Next, we
adaptively subdivide $\mtx{A}$ into disjoint blocks so that either $\omega_j r_k
\leq z$ or $\omega_j r_k > z$ for all $\omega_j$ and all $r_k$ in each block.
This leaves only a few small blocks with $\omega_j r_k \approx z$ whose entries
can be directly computed, and which can be directly applied.
Figure~\ref{fig:subdivide} shows a Hankel transform matrix $\mtx{A}$ divided
into local and asymptotic entries along the curve $\omega r = z$, as well as the
corresponding adaptive subdivision of the matrix into blocks which can be
rapidly applied. Following the subdivision step, all that remains is to apply
each of the disjoint blocks of $\mtx{A}$ to $\vct{f}$ using the corresponding
fast method.

%%% Local Variables: %% mode: latex %% TeX-master: "../main" %% End:

\section{Bessel function approximations} \label{sec:approx}
We now describe local and asymptotic expansions of the Bessel function
$J_\nu(\omega r)$, and provide error analysis by which one can select the number
of terms needed in each expansion to assure $\epsilon$ accuracy in both regimes.

\subsection{The Wimp expansion}\label{sec:local}

Near the origin, $J_\nu(z)$ is a smooth and essentially non-oscillatory function
of $z$. As a result, $J_\nu(xy)$ is a numerically low-rank function of all
sufficiently small inputs $x$ and $y$. Fortuitously, one such low-rank expansion
--- which we refer to as the \textit{Wimp expansion} --- is available in closed
form for integer~$\nu$~\cite{wimp1962polynomial}. In the case that $\nu$ is
even, we have
\begin{equation}
    \begin{aligned}
        J_\nu(xy) 
        &= \sum_{\ell=0}^\infty C_\ell(x) \, T_{2\ell}(y) \\
        C_\ell(x) 
        &= \delta_\ell \, J_{\frac{\nu}{2} + \ell}(x) \, J_{\frac{\nu}{2} - \ell}(x) \\
        \delta_\ell 
        &= \begin{cases}
            1 & \ell=0 \\
            2 & \text{otherwise}
        \end{cases}
    \end{aligned}
\end{equation}
for all $\abs{y} \leq 1$. A similar expansion exists for $\nu$
odd~\cite[2.23]{wimp1962polynomial}.

In order to employ the Wimp expansion to compute local terms within the Hankel
transform, we must determine the number of terms $L$ needed to construct an
$\epsilon$-accurate approximation to $J_\nu(\omega r)$ on a given rectangle
$(\omega, r) \in [0, \Omega] \times [0, R]$. The following lemma provides a
bound on the induced truncation error in the Wimp expansion as a function of the
order $\nu$, the space-frequency product $\Omega R$, and the number of retained
terms $L$.

\begin{lemma} \label{lem:wimp} Truncating the Wimp expansion after $L$ terms
    gives
    \begin{align} \label{eq:loc-error}
        \abs{J_\nu(\omega r) - \sum_{\ell=0}^L C_\ell(\omega R) T_{2\ell}\left( \frac{r}{R} \right)} 
        \leq \frac{2\exp\left\{ \frac{\nu}{2}(\beta - \gamma) + (L+1)(\beta + \gamma) \right\}}{1 - e^{\beta + \gamma}}
        =: \Bloc(\Omega R)
    \end{align}
    for all $\omega \in [0, \Omega], r \in [0, R]$, where
    \begin{align}
        \psi(p) &:= \log p + \sqrt{1 - p^2} - \log\left( 1 + \sqrt{1 - p^2} \right) \\
        \beta &:= \psi\left( \frac{\Omega R}{2L + 2 + \nu} \right) \\
        \gamma &:= \begin{cases}
            \psi\left( \frac{\Omega R}{2L + 2 - \nu} \right) & L + 1 \geq \frac{\nu}{2} \\
            0 & \textnormal{otherwise}
        \end{cases}
    \end{align}
\end{lemma}
\begin{proof}
    For $\nu$ even, the truncation error after $L$ terms is bounded by
    \begin{align}
        \abs{\sum_{\ell=L+1}^\infty C_\ell(\omega R) T_{2\ell}\left( \frac{r}{R} \right)}
        &\leq 2\sum_{\ell=L+1}^\infty \abs{J_{\frac{\nu}{2} + \ell}\left(\frac{\omega R}{2}\right)} \abs{J_{\frac{\nu}{2} - \ell}\left(\frac{\omega R}{2}\right)}.
    \end{align}
    Define $p_\ell(\omega) := \omega R / (\nu + 2\ell)$. Then by Siegel's bound
    \cite[10.14.5]{olver2010nist} we have
    \begin{align}
        \abs{J_{\frac{\nu}{2} + \ell}\left(\frac{\omega R}{2}\right)}
        &= \abs{J_{\frac{\nu}{2} + \ell}\bigg(\Big(\frac{\nu}{2} + \ell\Big) p_\ell(\omega)\bigg)} \\
        &\leq \exp\left\{ \Big(\frac{\nu}{2} + \ell\Big) \psi\big(p_\ell(\omega)\big)\right\} \\
        &\leq \exp\left\{ \Big(\frac{\nu}{2} + \ell\Big) \beta \right\},
    \end{align}
    where the last inequality follows from the fact that $\psi$ is an increasing
    function on $(0,1)$, and thus $\psi\big(p_\ell(\omega)\big) \leq \beta < 0$
    for all $\ell \geq L+1$ and all $\omega \in [0, \Omega]$. 
    
    If $L+1 \geq \frac{\nu}{2}$, we define $q_\ell(\omega) := \omega R / (2\ell
    - \nu)$ and apply Siegel's bound again to obtain
    \begin{align}
        \abs{J_{\frac{\nu}{2} - \ell}\left(\frac{\omega R}{2}\right)}
        &= \abs{J_{\ell - \frac{\nu}{2}}\bigg(\Big(\ell - \frac{\nu}{2}\Big) q_\ell(\omega)\bigg)} 
        \leq \exp\left\{ \Big(\ell - \frac{\nu}{2}\Big) \gamma \right\}.
    \end{align}
    If $L+1 < \frac{\nu}{2}$, Siegel's bound does not apply and we use instead
    the simple bound $\abs{J_{\frac{\nu}{2} - \ell}\left(\frac{\omega
    R}{2}\right)} \leq 1$, which is equivalent to taking $\gamma = 0$. 

    All that remains is to apply a geometric series argument
    \begin{align}
        \abs{J_\nu(\omega r) - \sum_{\ell=0}^L C_\ell(\omega R) T_{2\ell}\left( \frac{r}{R} \right)}
        &\leq 2\sum_{\ell=L+1}^\infty \exp\left\{ \Big(\frac{\nu}{2} + \ell\Big) \beta + \Big(\ell - \frac{\nu}{2}\Big) \gamma \right\} \\
        &= 2\exp\left\{\frac{\nu}{2}(\beta - \gamma)\right\} \sum_{\ell=L+1}^\infty \left( e^{\beta + \gamma} \right)^\ell \\
        &= \frac{2\exp\left\{ \frac{\nu}{2}(\beta - \gamma) + (L+1)(\beta + \gamma) \right\}}{1 - e^{\beta + \gamma}}
    \end{align}
    A similar calculation can be carried out for $\nu$ odd.
\end{proof}

Lemma~\ref{lem:wimp} is rather opaque regarding the impact of the various
parameters on the error because we have not utilized any simplifying bounds on
the function $\psi$, as done in~\cite[Lemma 1]{rangan2020factorization} for
large $\nu$. However, our analysis takes into account the decay in both
$J_{\frac{\nu}{2}+\ell}$ \textit{and} $J_{\frac{\nu}{2}-\ell}$, thus remaining
relatively tight for small $\nu$. It is therefore well-suited to our purposes
because, given $z, L > 0$, it provides a bound $\Bloc(z)$ on the pointwise error
in approximating any block of the matrix $J_\nu(\omega_j r_k)$ for which $\omega
r \leq z$ using the $L$-term Wimp expansion. 

This expansion is highly beneficial from a computational perspective, as it
yields an analytical rank-$L$ approximation to any block of $\bm{A}$ for which
$\omega_j r_k$ is sufficiently small
\begin{equation}
  \mtx{A}(j_0:j_1, k_0:k_1)
  %\vct{f}(k_0:k_1)
  \approx \mtx{C}\mtx{T}^\top
  %\vct{f}(k_0:k_1)
\end{equation}
where $\mtx{C} \in \R^{(j_1-j_0+1) \times L}$ and $\mtx{T} \in \R^{(k_1-k_0+1)
  \times L}$ with entries
\begin{equation}
  \mtx{C}(j,\ell) = C_{\ell-1}(\omega_j r_{k_1}) \qquad  \text{and} \qquad 
  \mtx{T}(k,\ell) = T_{2\ell-2}\left(\frac{r_k}{r_{k_1}}\right).
\end{equation}
For a block of $\mtx{A}$ of size $m_b \times n_b$, the low-rank approximation
given by the Wimp expansion can be applied to a vector in $\bO\big(L(m_b +
n_b)\big)$ time by first applying $\bm{T}^\top$ then applying $\bm{C}$.

\subsection{Hankel's expansion}
\label{sec:asymptotic}

Away from the origin, $J_\nu(z)$ exhibits essentially sinusoidal oscillation
with period $2\pi$. This statement is made precise by Hankel's asymptotic
expansion, which states that for $z \to \infty$
\begin{align} \label{eq:asymptotic-expansion}
    J_\nu(x)
    \sim \sqrt{\frac{2}{\pi x}} \left( 
        \cos\left(x + \phi\right) \sum_{\ell=0}^{\infty} \frac{(-1)^\ell a_{2\ell}(\nu)}{x^{2\ell}}
        - \sin\left(x + \phi\right) \sum_{\ell=0}^{\infty} \frac{(-1)^\ell a_{2\ell+1}(\nu)}{x^{2\ell+1}}
        \right)
\end{align}
where $\phi := - \frac{(2\nu+1)\pi}{4}$ and 
\begin{align}
    a_\ell(\nu) := \frac{(4\nu^2 - 1)(4\nu^2 - 3)\dots(4\nu^2 -
  (2\ell-1)^2)}{\ell! \, 8^\ell}.
\end{align}
Rearranging this expansion, we obtain an expansion which can be evaluated using
two NUFFTs and diagonal scalings, and whose remainder is bounded by the size of
the first neglected terms~\cite[Section 7.3]{watson1922treatise}
\begin{multline} \label{eq:asy-error}
    \Bigg| J_\nu(\omega r)
    - \sqrt{\frac{2}{\pi}} \sum_{\ell=0}^{M-1} \left[  
        \frac{(-1)^\ell a_{2\ell}(\nu)}{\omega^{2\ell + \frac{1}{2}}} \Re\left(
          \frac{e^{i(\omega r + \phi)}}{r^{2\ell+\frac{1}{2}}}\right)  - \frac{(-1)^\ell
            a_{2\ell+1}(\nu)}{\omega^{2\ell+\frac{3}{2}}}
          \Im\left(\frac{e^{i(\omega r + \phi)}}{r^{2\ell+\frac{3}{2}}} \right) \right]
        \Bigg| \\
         \leq \sqrt{\frac{2}{\pi}} \left( \frac{\abs{a_{2M}(\nu)}}{(\omega r)^{2M+\frac{1}{2}}} + \frac{\abs{a_{2M+1}(\nu)}}{(\omega r)^{2M+\frac{3}{2}}} \right) =: \Basy(z)
\end{multline}

The computational advantage of this expansion is that the $2M$-term asymptotic
expansion of any block of $\bm{A}$ can be rapidly applied to a vector~$\vct{x}$
using $2M$ Type-III NUFFTs
\begin{multline}
    \mtx{A}(j_0:j_1, k_0:k_1) \vct{x} 
    \approx \sqrt{\frac{2}{\pi}} \sum_{\ell=0}^{M-1} (-1)^\ell \Bigg[ 
        a_{2\ell}(\nu) \mtx{D}_\omega^{-2\ell-\frac{1}{2}} \Re\Big(e^{i\phi}
                                                             \mtx{F}
                                                             \mtx{D}_r^{-2\ell-\frac{1}{2}}
                                                             \vct{x} \Big) \\
        - a_{2\ell+1}(\nu) \mtx{D}_\omega^{-2\ell-\frac{3}{2}}
                                                             \Im\Big(e^{i\phi}
                                                               \mtx{F}
                                                               \mtx{D}_r^{-2\ell-\frac{3}{2}}
                                                               \vct{x} \Big) 
    \Bigg]
\end{multline}
where $\mtx{F} \in \C^{(j_1 - j_0 + 1) \times (k_1 - k_0 + 1)}$ is the Type-III
nonuniform DFT matrix corresponding to frequencies $\omega_{j_0}, \dots,
\omega_{j_1}$ and points $r_{k_0}, \dots, r_{k_1}$, and the diagonal scaling
matrices are given by $\mtx{D}_\omega :=
\diag(\omega_{j_0},\dots,\omega_{j_1})$, and $\mtx{D}_r :=
\diag(r_{k_0},\dots,r_{k_1})$.

\subsection{Determining order of expansions and crossover point}
\label{sec:cutoff}

With these error bounds in hand, we precompute the parameters~$z_{\nu,
\epsilon}^M$ and~$L_{\nu, \epsilon}^M$ for tolerances $\epsilon = 10^{-4},
\dots, \allowbreak 10^{-15}$, orders $\nu = 1, \dots, 100$, and number of
asymptotic expansion terms $M = 1, \dots, 20$:
\begin{itemize}
    \item $z_{\nu, \epsilon}^M$ such that $M$-term Hankel expansion of
    $J_\nu(\omega r)$ is $\epsilon$-accurate $\forall \ \omega r > z_{\nu,
    \epsilon}^M$,
    \item $L_{\nu, \epsilon}^M$ such that $L_{\nu, \epsilon}^M$-term Wimp
    expansion of $J_\nu(\omega r)$ is $\epsilon$-accurate $\forall \ \omega r
    \leq z_{\nu, \epsilon}^M$.
\end{itemize}
First, the crossover points $z_{\nu, \epsilon}^M$ are computed using Newton's
method on the function $\xi(z) := \Basy(z) - \epsilon$. Then the number of local
expansion terms $L_{\nu, \epsilon}^M$ are taken to be the smallest integer such
that $\Bloc\left(z_{\nu, \epsilon}^M\right) < \epsilon$. These tables are
precomputed once when the library is installed, and even this precomputation
requires only a few seconds on a laptop.

With these tables stored, for any order $\nu$ we can look up a pair of
complementary local and asymptotic expansions with error everywhere bounded by
the requested tolerance $\epsilon$. The only remaining free parameter is the
number of asymptotic terms $M$. This parameter is selected based on various
numerical experiments which maximize speed by balancing the cost of the local,
asymptotic, and direct evaluations. In our implementation, we use the heuristic
\begin{equation}
  \label{eq:num-asy-terms}
    M = \min\left(\flr{1 + \frac{\nu}{5} - \frac{\log_{10}(\epsilon)}{4}}, 20\right).
\end{equation}

%%% Local Variables: % mode: latex % TeX-master: "../main" % End:

\section{The Nonuniform Fast Hankel Transform} \label{sec:methods}

We now describe our NUFHT algorithm in detail, emphasizing the process by which
$\mtx{A}$ is adaptively subdivided into blocks using the results of the above
error analysis.

\subsection{Subdividing the matrix into blocks by expansion}

Having established error bounds which allow us to automatically select the
number of asymptotic terms~$M$, local terms~$L$, and crossover point~$z$ given a
tolerance~$\epsilon$ and order~$\nu$, we subdivide the matrix $\mtx{A}$ into
three sets of blocks, each of which can be efficiently applied to a vector as
described above:
\begin{itemize}
    \item Local blocks $\mathscr{L} = \big\{ \mtx{A}(j_0:j_1, k_0:k_1) \ | \
    \omega_j r_k \leq z \ \forall \ j_0 \leq j \leq j_1, \ k_0 \leq k \leq k_1
    \big\}$
    \item Asymptotic blocks $\mathscr{A} = \big\{ \mtx{A}(j_0:j_1, k_0:k_1) \ |
    \ \omega_j r_k > z \ \forall \ j_0 \leq j \leq j_1, \ k_0 \leq k \leq k_1
    \big\}$
    \item Direct blocks $\mathscr{D}$ which are small enough that no fast
    expansion is needed
\end{itemize}

In order to determine a subdivision of $\mtx{A}$ into blocks of these three
types, we initialize a set of \textit{mixed} blocks $\mathscr{M} = \{(1:m,
1:n)\}$, each of which contains a mix of local and asymptotic entries. We then
chose an index pair $(j,k)$ such that $\omega_j r_k \approx z$. This index
subdivides the block into four new sub-blocks with $(j,k)$ at the center, so
that the upper left block can be applied using the local expansion and is
appended to $\mathscr{L}$, and the lower right block using the asymptotic
expansion and is appended to $\mathscr{A}$. 

The remaining lower left and upper right blocks each still contain a mix of
local and asymptotic entries. If they are of sufficiently small size $m_b \times
n_b$ with $m_b n_b < \texttt{min\_size}$ --- a user-defined parameter which is
taken to be 1024 by default --- they can be evaluated directly and are appended
to $\mathscr{D}$. Otherwise they are appended back to $\mathscr{M}$, and we
continue the subdivision process recursively. 

This method yields a valid partition for any choice of $(j,k)$, but for
efficiency these indices are chosen to maximize the number of matrix entries
which can be applied using a fast expansion, i.e. the sizes of the upper left
and lower right blocks. This is done by solving the following constrained
optimization problem
\begin{align}
    (j,k) 
    &\ = \textsc{SplitIndices}(r_1,\dots,r_n, \omega_1,\dots,\omega_m, z) \\
    &:= \left\{
        \begin{array}{r@{\quad } l}
        \displaystyle\argmax_{j,k\in\Z} & (j-j_0)(k_1-k) + (j_1-j)(k-k_0)   \\
        \text{subject to} & j_0 \leq j \leq j_1 \\ 
        & k_0 \leq k \leq k_1 \\ 
        & \omega_j r_k \leq z
        \end{array}
    \right. \label{eq:subdiv-optim}
\end{align}
This problem can be solved exactly in $\bO(j_1-j_0 + k_1-k_0)$ time. However,
computing the exact optimal splitting indices for every box gives a negligible
speedup to the overarching Hankel transform compared to a simpler, quasi-optimal
scheme. In practice it is sufficient to choose a small number of equispaced
indices $j \in \{j_0,\dots,j_1\}$, compute the corresponding $k = \argmax \{k \
| \ r_k \leq \frac{z}{\omega_j}\}$ for each $j$, and choose $(j,k)$ as the pair
which minimizes the objective function of (\ref{eq:subdiv-optim}) among this
small collection.

\begin{algorithm2e}[t]
    \caption{Block subdivision of Hankel transform
    matrix}\label{alg:subdivision}
    \SetKwFunction{Subdivide}{Subdivide}
\Function{\Subdivide{$\bm{r}, \bm{\omega}, z, \texttt{\upshape min\_size}$}}{
    $\mathscr{L} = \mathscr{A} = \mathscr{D} = \emptyset$ \\
    $\mathscr{M} = \{(1:m, 1:n)\}$ \\
    \While{$\mathscr{M} \neq \emptyset$}{
        Pop an element $(j_0 : j_1, k_0 : k_1)$ from $\mathscr{M}$ \\
        $(j,k) = \textsc{SplitIndices}(r_{j_0},\dots,r_{j_1}, \omega_{k_0},\dots,\omega_{k_1}, z)$ \\
        Append $(j_0:j, k_0:k)$ to $\mathscr{L}$ \\
        Append $(j+1:j_1, k+1:k_1)$ to $\mathscr{A}$ \\
        Append $(j_0:j, k+1:k_1)$ to $\mathscr{M}$ \textbf{if} $(j - j_0 + 1)(k_1 - k) > \texttt{\upshape min\_size}$ \textbf{else} $\mathscr{D}$ \\
        Append $(j+1:j_1, k_0:k)$ to $\mathscr{M}$ \textbf{if} $(j_1 - j)(k_1 - k + 1) > \texttt{\upshape min\_size}$ \textbf{else} $\mathscr{D}$
    }
    \Return{$(\mathscr{L}, \mathscr{A}, \mathscr{D})$}
}
\end{algorithm2e}

\begin{algorithm2e}[t]
    \caption{Nonuniform fast Hankel transform}\label{alg:nufht}
    \SetKwFunction{NUFHT}{NUFHT}
\Function{\NUFHT{$\nu, \epsilon, \bm{r}, \bm{c}, \bm{\omega}$}}{
    $\bm{g} = \bm{0}$ \\
    Choose $M$ using (\ref{eq:num-asy-terms}) \\
    Look up $L = L_{\nu,\epsilon}^{M}$ and $z = z_{\nu,\epsilon}^{M}$ from pre-computed tables \\
    Choose \texttt{\upshape min\_size} from numerical experiments (e.g. default 1024)
    $(\mathscr{L}, \mathscr{A}, \mathscr{D}) = \textsc{Subdivide}(r_1,\dots,r_n, \omega_1,\dots,\omega_m, z, \texttt{\upshape min\_size})$ \\
    \For{$\mathscr{B} \in (\mathscr{L}, \mathscr{A}, \mathscr{D})$}{
        \For{$(j_0:j_1, k_0:k_1) \in \mathscr{B}$}{
            $\bm{g}(j_0:j_1) \pluseq \bm{A}(j_0:j_1, k_0:k_1) \bm{c}(k_0:k_1)$ using corresponding expansion \\
        }
    }
    \Return{$\bm{g}$}
}
\end{algorithm2e}

\subsection{Complexity analysis} \label{sec:complexity}

We now analyze the computational complexity of the proposed approach. In order
to do so, we must first comment on the complexity of the NUFFT, which is an
important subroutine in our method. Most analysis-based NUFFT codes ---
including the \texttt{FINUFFT} library \cite{barnett2019parallel} which we use
in our NUFHT implementation --- consist of three steps. First, delta masses
centered at each non-uniform point are convolved with a \textit{spreading
function} which smears them onto a fine $N$-point uniform grid. Then, a standard
equispaced FFT is computed on the fine grid. Finally, a diagonal de-convolution
with the Fourier transform of the spreading function is applied to reverse the
effect of the original smearing. For a more complete description of this NUFFT
method, see \cite{dutt1993fast,greengard2004accelerating,barnett2019parallel}.
For $n$ points $r_k$ and $m$ frequencies $\omega_j$, spreading the input points
to a finer grid is~$\bO(n)$, the FFT on the finer grid is $\bO(N\log N)$, and
the global deconvolution at the output frequencies is~$\bO(m)$. For the Type-III
NUFFT, the size $N$ of the fine grid typically scales linearly with the
space-frequency product $p := (\omega_m - \omega_1)(r_n - r_1)$
\cite{barnett2019parallel, greengard2004accelerating}. Therefore the total cost
of the NUFFT is $\bO(n + m + p\log p)$. Applying this fact in each asymptotic
block in the Hankel transform matrix, and adding the cost of applying local and
direct blocks, we can now analyze the complexity of the entire NUFHT method.

\begin{theorem} \label{thm:complexity} Take $\omega_1 < \dots < \omega_m \in
    [0,\infty)$ and $r_1 < \dots < r_n \in [0,\infty)$ and define the
    space-frequency product $p := (\omega_m - \omega_1)(r_n - r_1)$. Then the
    complexity of computing the NUFHT of order $\nu$ to tolerance $\epsilon$
    using Algorithm \ref{alg:nufht} is 
    $$\bO\Big((L + M)(m + n) \log \min(n,m) + Mp\log p\Big),$$ where $L$ and $M$
    are the number of local and asymptotic terms respectively chosen according
    to $\nu$ and $\epsilon$.
\end{theorem}

\begin{proof}
    For notational clarity we suppress the dependence of $z_{\nu, \epsilon}^M$
    on its parameters and simply denote it as $z$. If $\omega_j r_k \leq z$ for
    all $j=1,\dots,n$ and $k=1,\dots,m$ then only the $L$-term low-rank local
    expansion is used, which can be applied in $\bO(L(m + n))$ time. If instead
    $\omega_j r_k > z$ everywhere, then only the $M$-term asymptotic expansion
    is used, which can be applied using the Type-III NUFFT in $\bO(M(m + n +
    p\log p))$ complexity.

    Otherwise consider the case where $\mtx{A}$ contains both local and
    asymptotic entries. First, note that the number of levels $N_{\text{level}}$
    scales like $\bO(\log\min(n,m))$. The cost of determining the splitting
    indices $(j,k)$ for each box $\mtx{A}(j_0:j_1,k_0:k_1)$
    is~${\bO(j_1-j_0+k_1-k_0)}$, and thus the total cost of subdivision at each
    level is $\bO(m+n)$. Therefore the total cost of subdividing $\mtx{A}$ is
    $\bO((m+n)\log\min(n,m))$.
    
    Now, without loss of generality, assume~$\omega_1 \leq z/r_n < \omega_2$
     and~$r_1 \leq z/\omega_m < r_2$. If this were not the case, we would have
     blocks which can be evaluated using a single expansion as described above
     without affecting the complexity. After step~$\ell$ of subdividing every
     mixed block, we obtain $2^{\ell}$ new mixed blocks, $2^{\ell-1}$ new local
     blocks, and $2^{\ell-1}$ new asymptotic blocks. Let the local blocks be of
     size $m_{\ell,b}^{(\text{loc})} \times n_{\ell,b}^{(\text{loc})}$ for $b =
     1,\dots,2^{\ell-1}$. Then,
    \begin{equation}
      \sum_{b=1}^{2^{\ell-1}} m_{\ell,b}^{(\text{loc})} \leq m, \qquad
      \text{and} \qquad
      \sum_{b=1}^{2^{\ell-1}} n_{\ell,b}^{(\text{loc})} \leq n.
    \end{equation}
    An analogous fact holds for the asymptotic blocks.
    
    Therefore, the total cost of local evaluation is 
    \begin{equation}
      \begin{aligned}
        \sum_{\ell=1}^{N_{\text{level}}} \sum_{b=1}^{2^{\ell-1}} \bO\left(L\left(m_{\ell,b}^{(\text{loc})} + n_{\ell,b}^{(\text{loc})}\right)\right)
        &= \sum_{\ell=1}^{N_{\text{level}}} \bO(L(m + n)) \\
        &= \bO\big(L (m + n) \log \min(n,m)\big).
      \end{aligned}
    \end{equation}
    Let $p_{\ell,b}$ be the space-frequency product of box $b$ at level $\ell$.
    The total space frequency product $p$ is the area of the rectangle $R :=
    [\omega_1, \omega_m] \times [r_1, r_n]$, and all asymptotic boxes occupy
    disjoint sub-rectangles of $R$. Therefore the sum of their areas is bounded
    by the area of $R$, so that $$\sum_{\ell=1}^{N_{\text{level}}}
    \sum_{b=1}^{2^{\ell-1}} p_{\ell,b} \leq p.$$ Then by H\"older's inequality
    we obtain 
    \begin{align} 
        \sum_{\ell=1}^{N_{\text{level}}} \sum_{b=1}^{2^{\ell-1}} p_{\ell,b} \log p_{\ell,b}
        \leq \left( \sum_{\ell=1}^{N_{\text{level}}} \sum_{b=1}^{2^{\ell-1}} p_{\ell,b} \right) \left(\max_{\ell,b} \log p_{\ell,b} \right) 
        \leq p \log p.
    \end{align}
    The total cost of asymptotic evaluation via the Type-III NUFFT is therefore
    \begin{multline}
      \sum_{\ell=1}^{N_{\text{level}}} \sum_{b=1}^{2^{\ell-1}}
        \bO\left(M\left(m_{\ell,b}^{(\text{asy})} + n_{\ell,b}^{(\text{asy})} +
            p_{\ell,b}^{(\text{asy})} \log p_{\ell,b}^{(\text{asy})}
          \right)\right) \\
        \begin{aligned}
        &= \sum_{\ell=1}^{N_{\text{level}}} \bO(M(m + n)) +
        \sum_{\ell=1}^{N_{\text{level}}} \sum_{b=1}^{2^{\ell-1}}
        \bO\Big(M\Big(p_{\ell,b}^{(\text{asy})} \log
          p_{\ell,b}^{(\text{asy})}\Big)\Big)\\
        &= \bO\big(M (m + n) \log \min(n,m) + Mp\log p\big).
        \end{aligned}
      \end{multline}
    We subdivide until all direct blocks are all of size $m_b \times n_b$ with
    $m_bn_b = \bO(1)$. Thus the cost of computing the dense matvec with each
    direct block is $\bO(1)$, and the number of direct blocks is $\bO(m + n)$.
    Therefore the total direct evaluation cost is $\bO(m + n)$. Summing the cost
    of matrix subdivision, as well as local, asymptotic, and direct evaluation
    gives the result.
\end{proof}

In typical applications the maximum point $r_n$ is fixed by, for example, the
support of the function~$f$ whose Fourier transform is desired, and the maximum
frequency $\omega_m$ at which the transform is computed grows linearly with~$n$.
The following corollary studies this common scenario, which includes
Schl\"omilch expansions and Fourier-Bessel series. For notational conciseness,
we consider the number of terms~$L$ and~$M$ in each expansion as constants here.
\begin{corollary}
  \label{cor:complexity} Take $\omega_1 < \dots < \omega_n \in [0,\infty)$ and
    $r_1 < \dots < r_n \in [0,\infty)$ such that the space-frequency product $p
    = \bO(n)$. Then the complexity of computing the NUFHT using Algorithm
    \ref{alg:nufht} is $\bO(n\log n)$.
\end{corollary}

\begin{remark}
  There exist butterfly factorization-based NUFFT methods that could be used to
  remove the dependence on the space-frequency product $p$ in
  Theorem~\ref{thm:complexity} using linear algebraic
  approximations~\cite{pang2020interpolative}. However, we find that the
  asymptotic dependence on $p$ is generally seen only in pathological cases, and
  thus choose to avoid the precomputations associated with butterfly methods.

  % For local blocks away from the origin, one can prove that the
  % $\epsilon$-rank of blocks of $\mtx{A}$ is actually significantly lower than
  % the $L$ used in the Wimp expansion due to the complementary low-rank or
  % butterfly property of the Hankel transform~\cite{oneil2010algorithm}.
  % However, the evaluation of local blocks is neither an asymptotic nor
  % practical bottleneck, so we do not pursue an optimal-rank expansion in this
  % regime.
\end{remark}

%%% Local Variables: %% mode: latex %% TeX-master: "../main" %% End:

\section{Numerical experiments} \label{sec:results}
In the following section, we perform a number of numerical experiments to
validate the accuracy and complexity of our method. We close with two
applications from Fourier analysis and numerical PDEs.

\subsection{Comparison to direct evaluation}

We start by empirically verifying the error analysis in Sections~\ref{sec:local}
and \ref{sec:asymptotic}, and the asymptotic scaling analysis in
Section~\ref{sec:complexity} by comparing to direct evaluation of the Hankel
transform.

\subsubsection{Asymptotic scaling}
\begin{figure}[t]
  \centering
  \newcommand\twa{0.29cm} \newcommand\tw{0.43cm}
  \begin{subfigure}[b]{0.32\textwidth}
    \begin{tikzpicture}
        \draw (0, 0) node[inner sep=0] {\includegraphics[height=0.92\textwidth,
        trim={0.4cm 0 0.4cm 0}, clip]{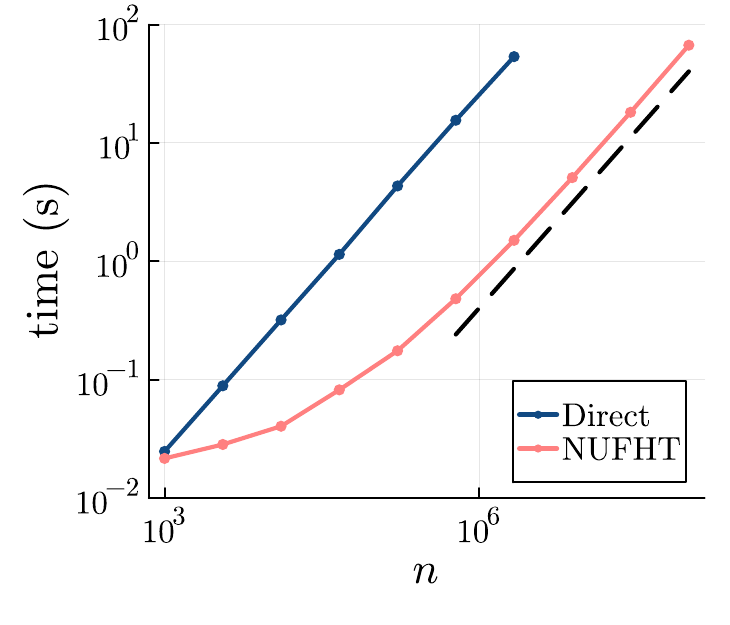}}; \draw (1.5,
        0.3) node {\small $\bO(n)$};
    \end{tikzpicture}
  \end{subfigure}
  \hfill
  \begin{subfigure}[b]{0.32\textwidth}
    \begin{tikzpicture}
        \draw (0, 0) node[inner sep=0] {\includegraphics[height=0.92\textwidth,
        trim={1.3cm 0 0.4cm 0}, clip]{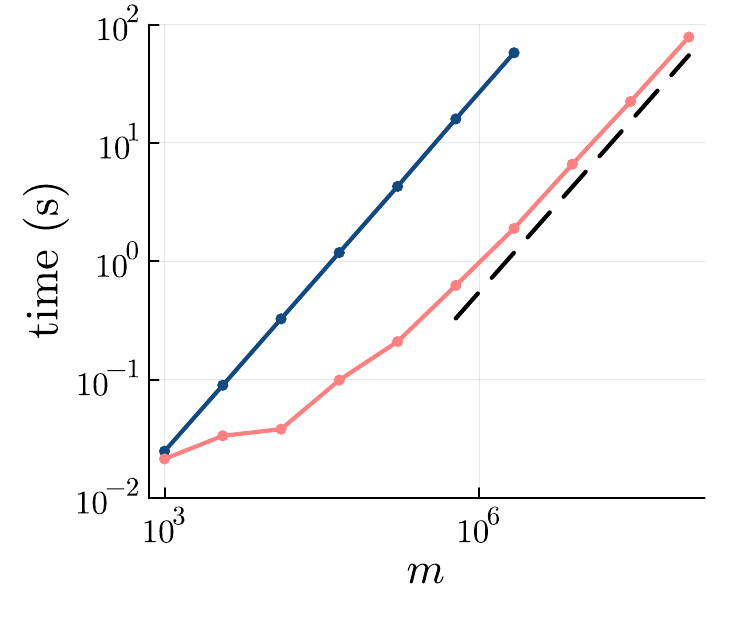}}; \draw (1.3,
        0.3) node {\small $\bO(m)$};
    \end{tikzpicture}
  \end{subfigure}
  % \hfill 
  \begin{subfigure}[b]{0.32\textwidth}
    \begin{tikzpicture}
        \draw (0, 0) node[inner sep=0] {\includegraphics[height=0.92\textwidth,
        trim={1.3cm 0 0.4cm 0}, clip]{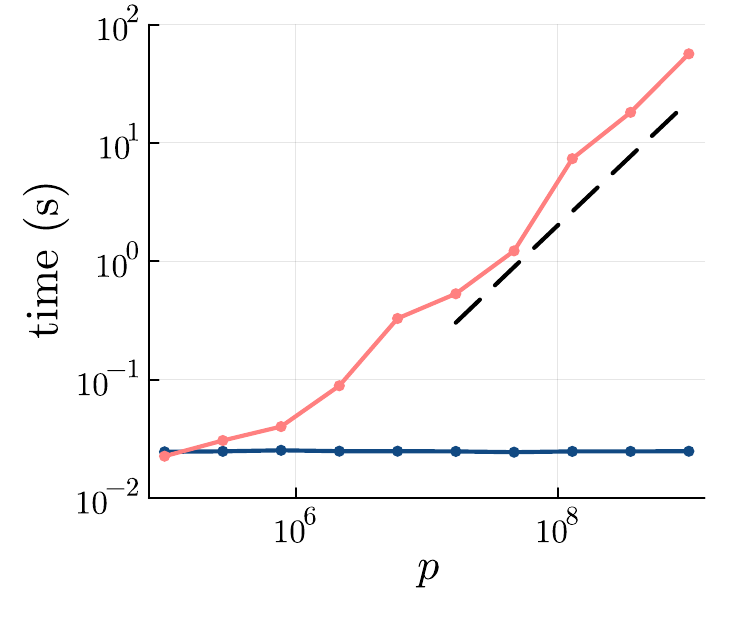}}; \draw (1.45,
        0.1) node {\small $\bO(p\log p)$};
    \end{tikzpicture}
  \end{subfigure}
  \caption{Scaling with $n$, $m$, and $p$ respectively, with the other variables
  held constant.}
  \label{fig:nmp-scaling}
\end{figure}
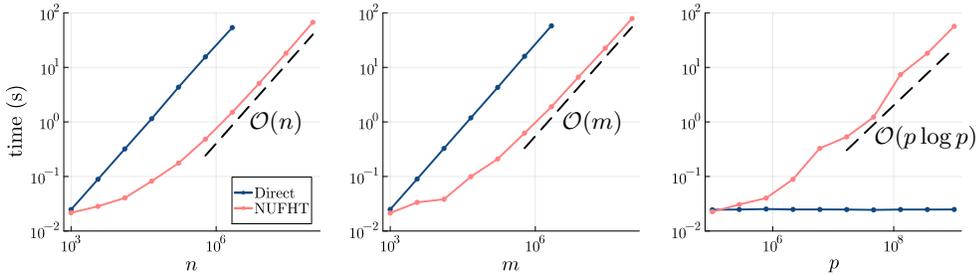

In order to study the impact of each of the relevant parameters in the scaling
analysis of Theorem~\ref{thm:complexity} independently, we take $n$ equispaced
points $r_k$ in the interval $[0,\sqrt{10^5}]$ and $m$ equispaced frequencies
$\omega_j$ in the interval $[0,p/\sqrt{10^5}]$. First, we fix $m=10^3$ and
$p=10^5$ while increasing $n$. Then, we fix~$n=10^3$ and $p=10^5$, this time
increasing $m$. Finally, we fix both $n = m = 10^3$ while increasing $p$.
Figure~\ref{fig:nmp-scaling} shows the CPU time for the NUFHT as well as for
direct summation in each of these scenarios. We observe the linear or
quasilinear scaling expected from Theorem~\ref{thm:complexity} with each of $n,
m,$ and $p$. Note in particular that the NUFHT scales with $p$ while direct
summation does not. Therefore, if a DHT is desired with relatively few points
with a very large space-frequency product, direct summation may give superior
performance, although such circumstances are rare in practice.

Next, we study the more typical scenario where the space-frequency product $p$
grows linearly with $n$, as discussed in Corollary~\ref{cor:complexity}. Here we
study two cases. First, we consider the Fourier-Bessel expansion where $\omega_j
= \mathrm{j}_{\nu, j}$ and $r_k = \mathrm{j}_{\nu, k}/\mathrm{j}_{\nu, n+1}$
with $n = m$. This is the direct analogue of the discrete Fourier transform as
the points and frequencies are the scaled roots of the basis, and the resulting
points and frequencies are quasi-equispaced for small to moderate $\nu$. 

We also consider the case of exponentially distributed points and frequencies
$\omega_j = r_j = 10^{\log_{10}(j) - \log_{10}(n)/2}$ with $n = m$. This is a
somewhat pathological worst case scenario for our algorithm, as the simple
calculation
\begin{equation}
  \sqrt{\frac{\Omega z}{R}} = \argmax_{\frac{z}{R} \leq \omega \leq \Omega} \ (\Omega - \omega)\left(R - \frac{z}{\omega}\right)
\end{equation}
shows that if we subdivide a block with space frequency product $\Omega R$ at a
point $(\omega, r)$ which lies on the curve $\omega r = z$, then the largest
possible space-frequency product~$p$ for the resulting lower right asymptotic
block is achieved by taking $\omega$ to be the mid-point of $[z/R, \Omega]$ on a
log scale. In other words, points and frequencies which are exponentially
distributed result in the highest possible space-frequency product $p$ for every
asymptotic block at every level. From Theorem~\ref{thm:complexity}, maximizing
$p$ drives the cost of the NUFHT. This distribution of points and frequencies is
also challenging because it leads to equally-sized square blocks at every level,
which guarantees that all blocks are subdivided the maximum number of times
before yielding sufficiently small direct blocks.

Figure~\ref{fig:both-scaling} shows the CPU time needed to evaluate the NUFHT in
the Fourier-Bessel and exponentially-distributed cases with $\nu=0$ and
$\epsilon=10^{-8}$. Both cases eventually demonstrate the expected $\bO(n\log
n)$ scaling. As a result of the challenges just discussed for the
exponentially-distributed case, its runtime is up to an order of magnitude
slower than the Fourier-Bessel series.

\subsubsection{Impact of the order and tolerance on runtime} 
\begin{figure}
  \centering
  \begin{subfigure}[b]{0.32\textwidth}
    \begin{tikzpicture} 
        \draw (0, 0) node[inner sep=0] {\includegraphics[height=0.92\textwidth,
        trim={0.4cm 0 0.4cm 0}, clip]{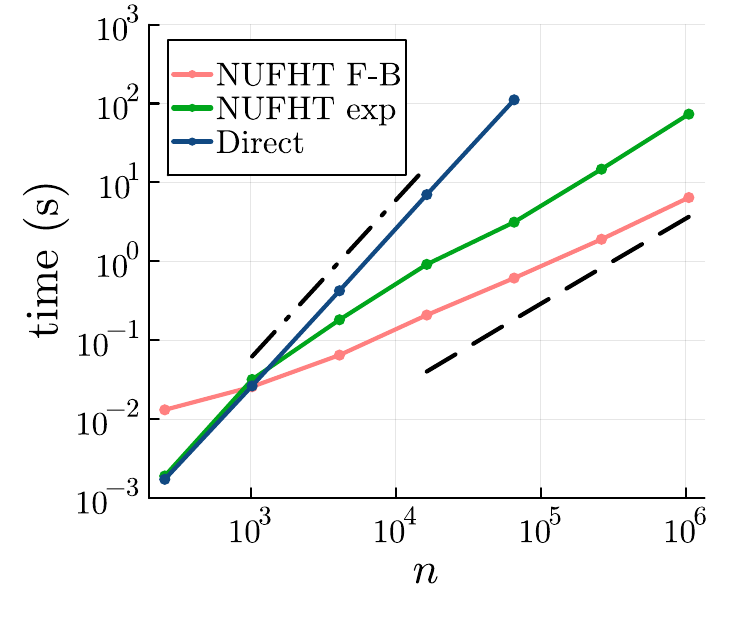}}; \draw (1.25,
        -0.5) node {\small $\bO(n\log n)$}; \draw (-0.7, 0.4) node {\small
        $\bO\big(n^2\big)$};
    \end{tikzpicture}
  \end{subfigure}
  \hfill
  \begin{subfigure}[b]{0.32\textwidth}
    \includegraphics[height=0.92\textwidth, trim={1.3cm 0 0.4cm 0}, clip]{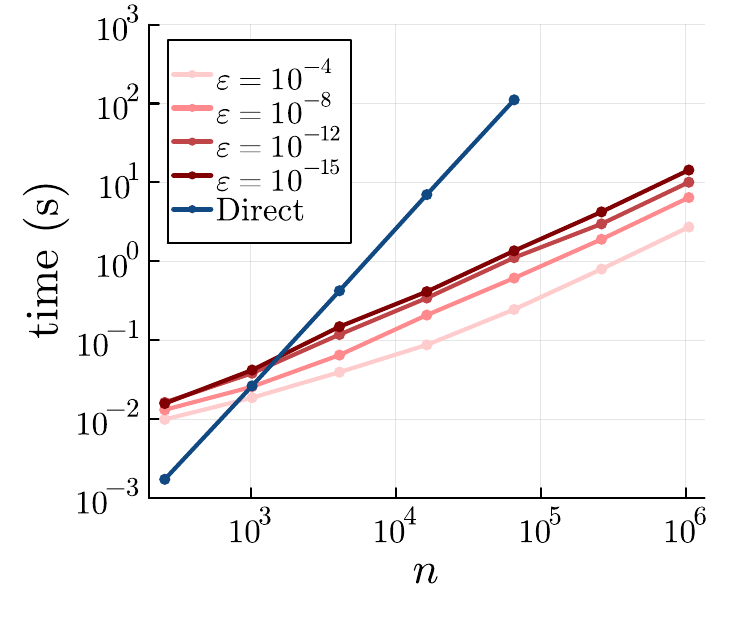}
  \end{subfigure} 
  % \hfill 
  \begin{subfigure}[b]{0.32\textwidth}
    \includegraphics[height=0.92\textwidth, trim={1.3cm 0 0.4cm 0}, clip]{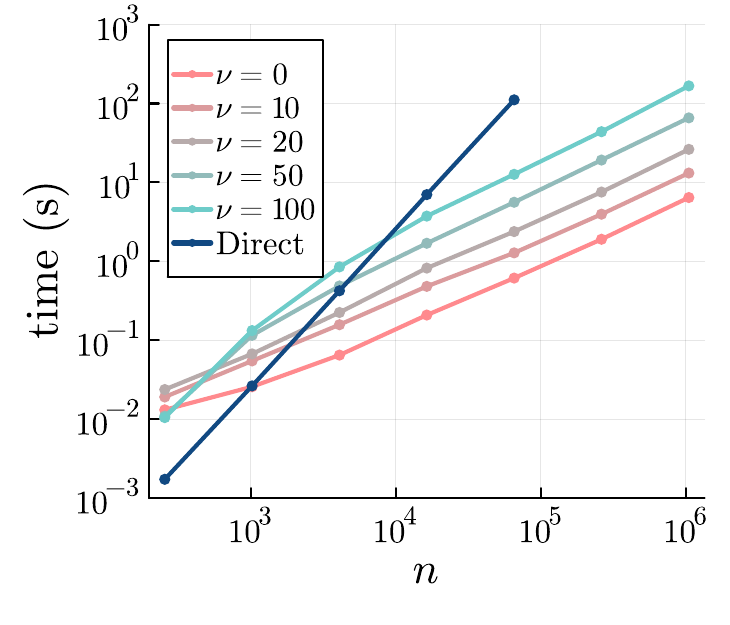}
  \end{subfigure}
  \caption{Scaling with $n$ for $p = \bO(n)$ test cases. In the first plot, we
  fix $\nu = 0, \epsilon = 10^{-8}$ and time the NUFHT for both the
  Fourier-Bessel and exponentially distributed cases. In the second and third
  plots, we consider the Fourier-Bessel series only, and fix one of the
  parameters $\nu = 0$ and $\epsilon = 10^{-8}$ while varying the other. The
  timings of direct summation and Fourier-Bessel series from the first plot are
  repeated in the other two plots for reference.}
  \label{fig:both-scaling}
\end{figure}

As the order $\nu$ increases or the tolerance $\epsilon$ decreases, the number
of necessary terms $L$ and $M$ in the local and asymptotic expansions,
respectively, both grow. From Theorem~\ref{thm:complexity}, we expect the
runtime to grow linearly with $L + M$. Figure~\ref{fig:both-scaling} shows the
runtime of our method for various $\epsilon$ with $\nu=0$ held constant, as well
as for multiple $\nu$ with $\epsilon=10^{-8}$ fixed. The $\bO(n\log n)$ scaling
of the algorithm is similar in all cases, while the prefactors vary; a transform
with $\epsilon = 10^{-15}$ is about an order of magnitude slower than using
$\epsilon = 10^{-4}$, and an order $\nu=100$ transform is almost two orders of
magnitude slower than the order $\nu=0$ equivalent.

\subsubsection{Approximation error}

\begin{figure}[t]
  \hspace*{0.25\textwidth}
  \includegraphics[width=0.4\textwidth]{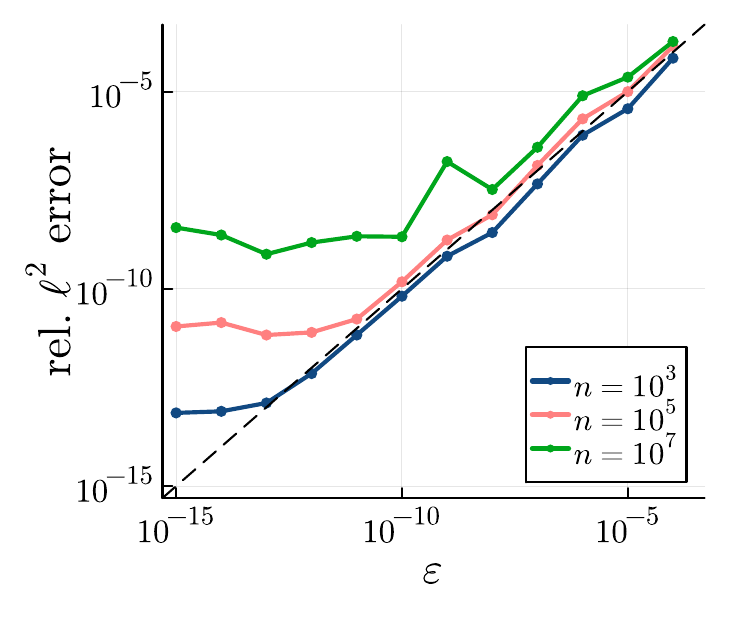}
  \caption{Relative 2-norm error $\norm[2]{\vct{g} - \vct{\tilde{g}}} /
  \norm[2]{\vct{g}}$ as a function of tolerance $\epsilon$ for a NUFHT of order
  $\nu=0$ for various $n$.}
  \label{fig:accuracy}
\end{figure}

Finally, we study the relative error in the output $\vct{g}$ as a function of
the desired tolerance $\epsilon$. To do this, we fix $n$ and form a sparse
vector $\vct{f} \in \R^n$ with 1000 nonzero entries whose indices are selected
at random and whose values are independent standard Gaussian. We evaluate the
Fourier-Bessel series using the NUFHT with the full vector~$\vct{f}$ as input,
and denote the output as $\vct{\tilde{g}}$. We then use direct summation on only
the nonzero entries to generate a reference result $\vct{g}$.
Figure~\ref{fig:accuracy} shows the 2-norm relative error $\norm[2]{\vct{g} -
\vct{\tilde{g}}} / \norm[2]{\vct{g}}$ between the NUFHT and the reference. For
small transforms with $n=10^3$, the relative error demonstrates excellent
agreement with the tolerance $\epsilon$ down to $\epsilon = 10^{-14}$ or so.
This suggests that the analysis used in Section~\ref{sec:approx} to determine
the necessary number of local and asymptotic terms is fairly tight. For larger
transforms, however, the error saturates, and regardless of the tolerance
$\epsilon$ our method gives at most 9 digits of accuracy for transforms of size
$n=10^7$. This is a well-known limitation of existing NUFFT methods, for which
the error generally scales like $n$ times machine precision \cite[Remark
9]{barnett2019parallel}.

\subsection{Computing Fourier transforms of radial functions}

For radial functions $f(\bm{r}) = f(\norm{\bm{r}})$ in $\R^d$, one can integrate
out the radial variables analytically, reducing the $d$-dimensional Fourier
integral to a single Hankel transform
\begin{align} \label{eq:radial-fourier}
    \hat{f}(\bm{\omega}) 
    = \int_{\R^d} f(\norm{\bm{r}}) \, e^{i\bm{\omega}^\top \bm{r}} \, \dif{\bm{r}}
    = \frac{(2\pi)^{\frac{d}{2}}}{\omega^{\frac{d}{2} - 1}} \int_0^\infty f(r)
  \, J_{\frac{d}{2} - 1}(\omega r) \, r^{\frac{d}{2}} \, \dif{r}.
\end{align}

We compare two methods of computing $\hat{f}$ for the indicator function of the
unit disk $f(r) = \ind{0 \leq r \leq 1}$ to absolute error $\epsilon = 10^{-12}$
at $n$ equispaced points $\omega_j \in [0, \omega_{\text{max}}]$. First, we use
a Gauss-Legendre quadrature rule on $[0,1]$ with nodes~$r_k$ and weights~$w_k$.
We utilize the NUFHT to compute the resulting sum
\begin{equation}
  \begin{aligned}
  \hat{f}(\omega) 
  &= 2\pi\int_0^1 f(r) \, J_0(\omega r) \, r \, \dif{r} \\
  &\approx 2\pi \sum_{k=1}^m w_k \, f(r_k) \, J_0(\omega r_k) \, r_k,
  \end{aligned}
\end{equation}
doubling the number of nodes~$m$ until the error in the computed integral is
less than~$\epsilon$. Second, we build a two-dimensional quadrature rule in
polar coordinates, using the same $m$-point Gauss-Legendre rule in $r$ and a
$t_k$-node trapezoidal rule in $\theta$ on each circle of radius $r_k$. We
double the number of trapezoidal nodes $t_k$ in each circle until the error in
the corresponding radial integral is less than $\epsilon$. We then utilize the
2D NUFFT to compute the resulting double sum
\begin{align}
  \hat{f}(\omega) 
  &= \frac{1}{4\pi^2} \int_0^{2\pi} \int_0^1 f(r) \, e^{-i\omega r\cos\theta} \,
    r \, \dif{r} \dif{\theta} \\
  &\approx \frac{1}{4\pi^2} \sum_{k=1}^{m} w_k \, r_k \, f(r_k) \,
    \frac{2\pi}{t_k}
    \sum_{s=1}^{t_k} \exp\left\{-i\omega r_k \, \cos\left(\frac{2\pi s}{t_k}\right)\right\}.
\end{align}

\begin{figure}
  \centering
  \begin{subfigure}[b]{0.38\textwidth}
    \includegraphics[height=\textwidth]{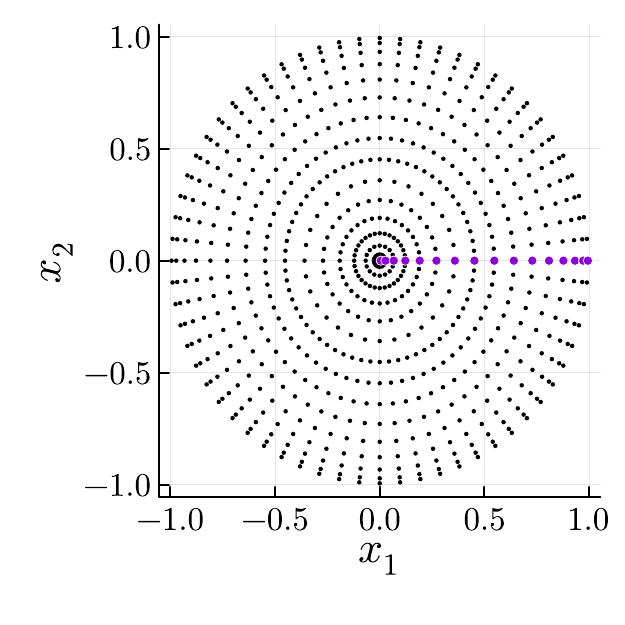}
  \end{subfigure}
  % \hfill
  \begin{subfigure}[b]{0.60\textwidth}
    \begin{tikzpicture} 
      \draw (0, 0) node[inner sep=0]
      {\includegraphics[height=0.63\textwidth]{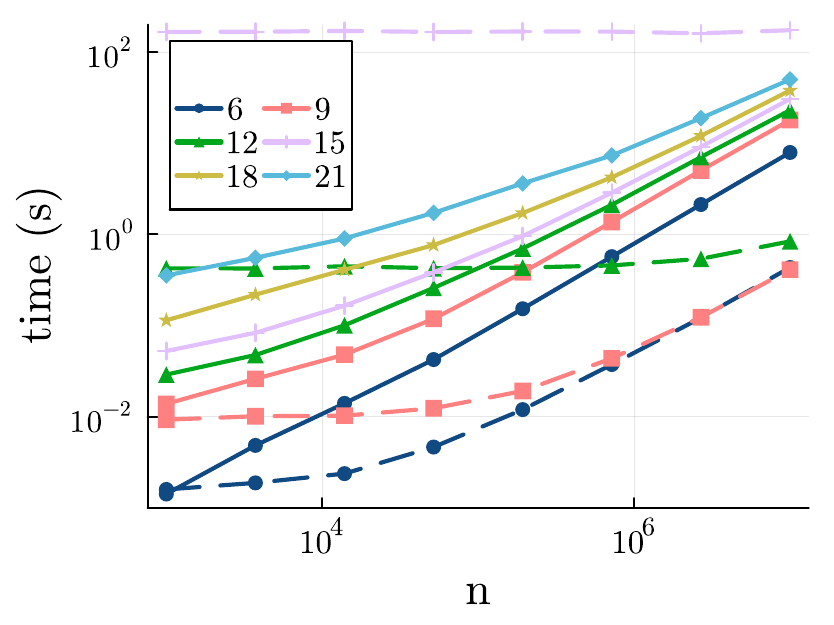}};
      \draw (-1.2, 1.92) node {\scriptsize $\log_2\omega_{\text{max}}$};
    \end{tikzpicture}
  \end{subfigure}
  \caption{Example two-dimensional quadrature nodes for the NUFFT, with
  one-dimensional radial Gauss-Legendre quadrature on $[0,1]$ for the NUFHT
  emphasized. Runtime comparison between NUFHT and 2D NUFFT for various choices
  of the maximum frequency $\omega_{\text{max}}$ and the number of  evaluation
  points $n$. Solid lines indicate the NUFHT, and the corresponding dashed lines
  indicate the 2D NUFFT.}
  \label{fig:fourier-test}
\end{figure}

If only low frequencies $\omega$ are desired, e.g. $\omega_{\text{max}} = 64$,
the integrands are only mildly oscillatory and few trapezoidal nodes are
required. In combination with the relative ease of amortizing costs in the
NUFFT, the two-dimensional transform is often faster than the NUFHT. However,
for larger $\omega_{\text{max}}$ the integrands become more oscillatory, and in
two dimensions $m = \bO(\omega_{\text{max}}^2)$ nodes are needed to resolve
these oscillations. Therefore the $\bO(m)$ spreading step in the NUFFT becomes
prohibitively expensive. However, by using radial symmetry to reduce to a
one-dimensional integral, the NUFHT requires only $\bO(\omega_{\text{max}})$
quadrature nodes, avoiding the curse of dimensionality.
Figure~\ref{fig:fourier-test} shows an example quadrature and runtimes for both
the NUFFT and NUFHT approaches. Note that for $\omega_{\text{max}} = 2^{15}$ the
2D NUFFT is orders of magnitude slower than the NUFHT for most $n$, and for even
larger $\omega_{\text{max}}$ the quadratic scaling of the 2D NUFFT with
frequency makes the computation intractable on a laptop, while the NUFHT's
linear scaling with frequency allows evaluation of the Fourier transform at
significantly higher frequencies at an only moderately increased cost.

\subsection{A Helmholtz solver using Fourier-Bessel expansions}

Finally, we demonstrate the application of the nonuniform Hankel transform to
solving partial differential equations on the disk using Fourier-Bessel
expansions. Consider the following inhomogeneous Helmholtz problem on the unit
disk~$D$
\begin{equation}
  \label{eq:helmholtz}
    \begin{aligned}
      (\Delta + \kappa^2) u(r,\theta) &= f(r,\theta), & \qquad &\text{for } r\in
                                                                 [0,1), \quad
                                                                 \theta \in
                                                                 [0,2\pi), \\
      u(1,\theta) &= 0 & \qquad &\text{for } \theta \in [0,2\pi).
    \end{aligned}
\end{equation}
Note that the functions $\psi_{j\ell}(r,\theta) := J_\ell(\mathrm{j}_{\ell,j} r)
e^{i\ell\theta}$ are the eigenfunctions of the Laplacian on the unit disk with
homogeneous Dirichlet boundary condition, so that
\begin{equation}
  \Delta \psi_{j\ell}(r,\theta)
  = \lambda_{j\ell} \psi_{j\ell}(r,\theta),
\end{equation}
where~$\lambda_{j\ell} = -
\mathrm{j}_{\ell,j}^2$~\cite{boyd2011comparing,watson1922treatise}. Therefore,
writing the forcing function~$f$ and solution~$u$ in terms of their respective
Fourier-Bessel expansions
\begin{equation}
  f(r,\theta) 
  = \sum_{\ell=-\infty}^\infty \sum_{j=1}^\infty \alpha_{j\ell} \, 
  J_\ell(\mathrm{j}_{\ell,j}r) \, e^{i\ell\theta}, \qquad
  u(r,\theta) 
  = \sum_{\ell=-\infty}^\infty \sum_{j=1}^\infty \beta_{j\ell} \,
  J_\ell(\mathrm{j}_{\ell,j}r)
  \, e^{i\ell\theta}
\end{equation} 
decouples~\eqref{eq:helmholtz} into a system of diagonal equations resulting in
an explicit formula for the coefficients~$\beta_{j\ell}$:
%\begin{align} \beta_{j\ell} (\Delta + \kappa^2) \psi_{j\ell}(r, \theta) =
%  \beta_{j\ell} (\lambda_{j\ell} + \kappa^2) \psi_{j\ell}(r, \theta) &=
%  \alpha_{j\ell} \psi_{j\ell}(r, \theta) \end{align} which implies
\begin{equation}
  \beta_{j\ell} 
  = \frac{\alpha_{j\ell}}{\lambda_{j\ell} + \kappa^2}.
\end{equation}
Due to the orthogonality of the Bessel functions~$J_\ell$, the Fourier-Bessel
coefficients of the forcing~$f$ can be computed as:
\begin{equation} \label{eq:FB-coef}
  \alpha_{j\ell} 
  = \frac{2}{J_{\ell+1}(\mathrm{j}_{\ell,j})^2} \int_0^{2\pi} \int_0^1
  f(r,\theta) \, J_\ell(\mathrm{j}_{\ell,j}r) \, e^{-i\ell\theta} \, r \dif{r} \dif{\theta},
\end{equation}
and the Fourier-Bessel expansion of the solution $u$ can then be written
explicitly
\begin{equation}
  u(r,\theta) = \sum_{\ell=-\infty}^\infty \sum_{j=1}^\infty
  \frac{\alpha_{j\ell}}{\lambda_{j\ell} + \kappa^2} \,
  J_\ell(\mathrm{j}_{\ell,j}r) \, e^{i\ell\theta}.
\end{equation}
By diagonalizing the Laplacian, this Fourier-Bessel solver thus provides a
direct analogue in the Dirichlet disk setting of spectral methods on a periodic
rectangle using bivariate Fourier expansions, and inherits many of the merits of
spectral methods. First, if~$f$ and all its derivatives go to zero at~$r=1$
and~$f$ is smooth in the interior of~$D$, then $\abs{\alpha_{j\ell}} \to 0$
exponentially fast in both $j$ and $\ell$ \cite{boyd2011comparing}. In addition,
solutions for arbitrary~$\kappa$ can be evaluated without additional
computations involving $f$, assuming that~$\kappa^2$ is not itself a Dirichlet
eigenvalue of the Laplacian on~$D$.

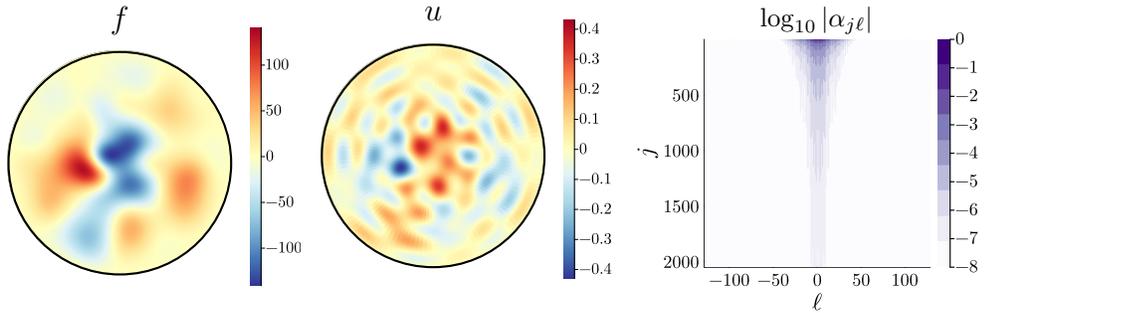
\begin{figure}[t]
  \centering
  \begin{subfigure}[T]{0.3\textwidth}
    \begin{tikzpicture}
      \draw (0, 0) node[inner sep=0] {\includegraphics[width=\textwidth,
      trim={2cm 1cm 0.1cm 0.2cm}, clip]{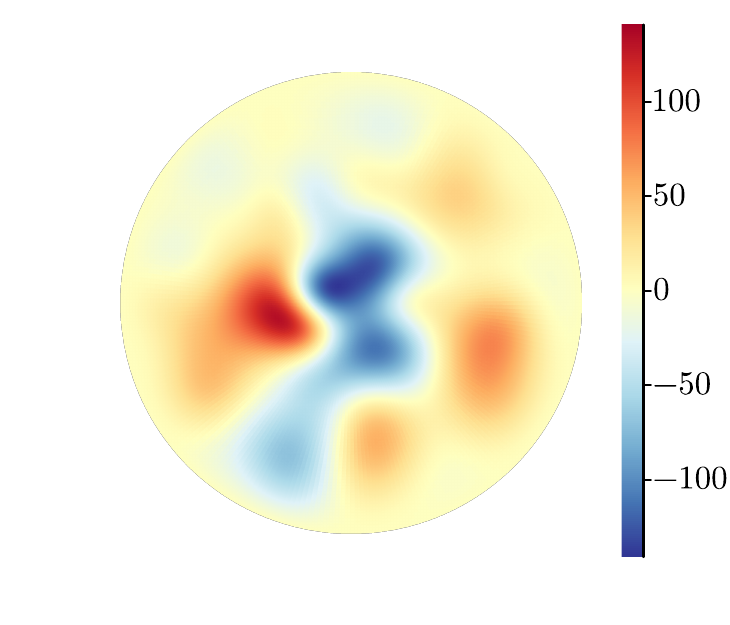}};
      \draw[fill=none, thick](-0.44, -0.1) circle (1.48) node {}; \draw[](-0.44,
      1.8) node {\large $f$};
    \end{tikzpicture}
  \end{subfigure}
  \hfill
  \begin{subfigure}[T]{0.3\textwidth}
    \begin{tikzpicture}
      \draw (0, 0) node[inner sep=0] {\includegraphics[width=\textwidth,
      trim={2cm 1cm 0.1cm 0.2cm}, clip]{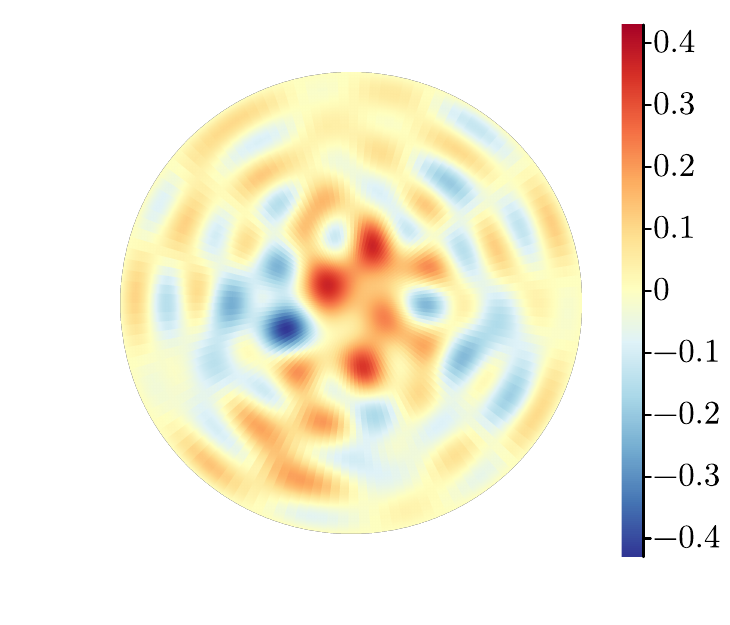}};
      \draw[fill=none, thick](-0.44, -0.1) circle (1.48) node {}; \draw[](-0.44,
      1.8) node {\large $u$};
    \end{tikzpicture}
  \end{subfigure} 
  \hfill 
  \begin{subfigure}[T]{0.36\textwidth}
    \begin{tikzpicture}
      \draw (0, 0) node[inner sep=0] {\includegraphics[width=\textwidth,
      trim={0.9cm 1cm 0.5cm 0.5cm}, clip]{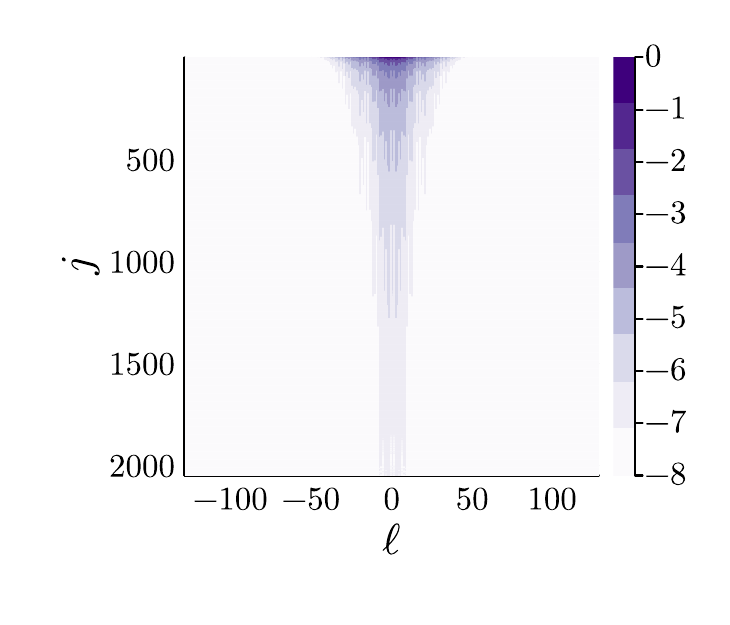}};
      \draw[](0.1, 2) node {$\log_{10}\abs{\alpha_{j\ell}}$};
    \end{tikzpicture}
  \end{subfigure}
  \caption{Forcing $f$, solution $u$, and log magnitude of Fourier-Bessel
  expansion coefficients $\alpha_{j\ell}$ for \eqref{eq:helmholtz} with $\kappa
  = 25$.}
  \label{fig:helmholtz}
\end{figure}

To compute the Fourier-Bessel coefficients $\alpha_{j\ell}$ of $f$
using~\eqref{eq:FB-coef}, we use an $m$-point Gauss-Legendre rule in $r$ and a
$t$-point trapezoidal rule in $\theta$. We iteratively double the number of
nodes in each rule until the relative norm difference in computed coefficients
between iterations is less than $\epsilon$ (controlling the discretization
error) and the relative norm of the coefficients appended in the last iteration
is less than $\epsilon$ (controlling the truncation error). Computing all
$\alpha_{j\ell}$ at each iteration requires $t$ NUFHTs of size $m$ and $m$ FFTs
of size $t$, resulting in $\bO(tm\log m + mt\log t)$ total complexity.
Figure~\ref{fig:helmholtz} shows an example random forcing $f$, the magnitude of
its Fourier-Bessel coefficients $\alpha_{j\ell}$, and the corresponding solution
$u$ to the Helmholtz equation~\eqref{eq:helmholtz} compute to relative precision
$\epsilon = 10^{-8}$.

This approach does, however, have two main limitations. First is that the
coefficients of $f$ decrease only algebraically in $j$ if $f$ has nonzero
derivatives at~$r=1$. More precisely, if $\Delta^q f(r) |_{r=1} = 0$ for all $0
\leq q \leq p-1$, then $\abs{\alpha_{j\ell}} \sim j^{-2p - \frac{1}{2}}$, with
exponential convergence only possible if $\Delta^q f(r) |_{r=1} = 0$ for all
integer $q$~\cite{boyd2011comparing}.  This is a fundamental property of the
Fourier-Bessel expansion, and does not depend on the numerical method used to
evaluate the Hankel transform. The second limitation is the increase in
computational cost of our NUFHT with the order $\ell$, as demonstrated in
Figure~\ref{fig:both-scaling}. As $\alpha_{j\ell}$ decrease spectrally in $\ell$
for smooth functions $f$, very large $\ell$ are not often needed. However, as in
any spectral method, functions with sharp features or discontinuous derivatives
will yield only algebraic decay in~$\ell$, requiring more Fourier bases. In such
cases the corresponding high order NUFHTs become intractable using the method
described here. 

%%% Local Variables: %% mode: latex %% TeX-master: "../main" %% End:

\section{Discussion} \label{sec:discussion}

In this manuscript we have presented a fast algorithm for computing discrete Hankel
transforms of moderate orders from $n$ nonuniform points to $m$ nonuniform
frequencies in $\bO\big((m+n)\log\min(n,m)\big)$ operations. The algorithm
relies on a careful space-frequency analysis of the Bessel function kernel,
judicious use of small-argument series expansions and large-argument asymptotic
expansions, as well as a small number of direct calculations. The algorithm
makes no assumptions on the distribution of points in space and frequency --- it
applies to the fully nonuniform case --- and can be used for Hankel transforms
of higher order with a modest increase in computational cost. More importantly,
the algorithm does not require any precomputation, in contrast to algorithms
based on butterfly factorizations of the Hankel transform matrix. Significant
speedups over the direct calculation have been demonstrated, as well as
asymptotic scaling of the computational complexity. An implementation of the
algorithm of this paper is available as an open-source Julia package at
\href{https://github.com/pbeckman/FastHankelTransform.jl}{\texttt{github.com/pbeckman/FastHankelTransform.jl}}.

In order to efficiently extend our algorithm to compute arbitrarily high-order
Hankel transforms which are needed for higher-order Fourier-Bessel expansions
and in various high-dimensional statistical
settings~\cite{lord1954a,nolan2013multivariate}, alternative expansions and
asymptotics of~$J_\nu$ need to be used or derived. This is the focus of ongoing
research.

%%% Local Variables: % mode: latex % TeX-master: "../main" % End:

\section*{Acknowledgments}
The authors would like to thank Alex Barnett for suggesting the use of the Wimp
expansion.

\section*{Competing interests}
The authors report no competing interests.

\bibliographystyle{siamplain}
\bibliography{refs}

\begin{thebibliography}{10}

\bibitem{alexander2012adaptive}
{\sc T.~S. Alexander}, {\em Adaptive signal processing: theory and applications}, Springer Science \& Business Media, 2012.

\bibitem{ali1999generalized}
{\sc I.~Ali and S.~Kalla}, {\em A generalized {H}ankel transform and its use for solving certain partial differential equations}, The ANZIAM Journal, 41 (1999), pp.~105--117.

\bibitem{alpert2002adaptive}
{\sc B.~Alpert, G.~Beylkin, D.~Gines, and L.~Vozovoi}, {\em Adaptive solution of partial differential equations in multiwavelet bases}, Journal of Computational Physics, 182 (2002), pp.~149--190.

\bibitem{askham2017adaptive}
{\sc T.~Askham and A.~J. Cerfon}, {\em An adaptive fast multipole accelerated poisson solver for complex geometries}, Journal of Computational Physics, 344 (2017), pp.~1--22.

\bibitem{barnett2019parallel}
{\sc A.~H. Barnett, J.~Magland, and L.~af~Klinteberg}, {\em A parallel nonuniform fast {F}ourier transform library based on an “exponential of semicircle" kernel}, SIAM Journal on Scientific Computing, 41 (2019), pp.~C479--C504.

\bibitem{bisseling1985fast}
{\sc R.~Bisseling and R.~Kosloff}, {\em The fast {H}ankel transform as a tool in the solution of the time dependent {S}chr{\"o}dinger equation}, Journal of Computational Physics, 59 (1985), pp.~136--151.

\bibitem{bondesson2019nonuniform}
{\sc D.~Bondesson, M.~J. Schneider, T.~Gaass, B.~K{\"u}hn, G.~Bauman, O.~Dietrich, and J.~Dinkel}, {\em Nonuniform {F}ourier-decomposition {MRI} for ventilation-and perfusion-weighted imaging of the lung}, Magnetic resonance in medicine, 82 (2019), pp.~1312--1321.

\bibitem{boyd2011comparing}
{\sc J.~P. Boyd and F.~Yu}, {\em Comparing seven spectral methods for interpolation and for solving the {P}oisson equation in a disk: {Z}ernike polynomials, {L}ogan--{S}hepp ridge polynomials, {C}hebyshev--{F}ourier series, cylindrical {R}obert functions, {B}essel--{F}ourier expansions, square-to-disk conformal mapping and radial basis functions}, Journal of Computational Physics, 230 (2011), pp.~1408--1438.

\bibitem{bronstein2002reconstruction}
{\sc M.~M. Bronstein, A.~M. Bronstein, M.~Zibulevsky, and H.~Azhari}, {\em Reconstruction in diffraction ultrasound tomography using nonuniform {FFT}}, IEEE transactions on medical imaging, 21 (2002), pp.~1395--1401.

\bibitem{brunol1977fourier}
{\sc J.~Brunol and P.~Chavel}, {\em Fourier transformation of rotationally invariant two-variable functions: {C}omputer implementation of {H}ankel transform}, Proceedings of the IEEE, 65 (1977), pp.~1089--1090.

\bibitem{cavanagh1979numerical}
{\sc E.~Cavanagh and B.~Cook}, {\em Numerical evaluation of {H}ankel transforms via {G}aussian-{L}aguerre polynomial expansions}, IEEE transactions on acoustics, speech, and signal processing, 27 (1979), pp.~361--366.

\bibitem{cree1993algorithms}
{\sc M.~Cree and P.~Bones}, {\em Algorithms to numerically evaluate the {H}ankel transform}, Computers \& Mathematics with Applications, 26 (1993), pp.~1--12.

\bibitem{dutt1993fast}
{\sc A.~Dutt and V.~Rokhlin}, {\em Fast {F}ourier transforms for nonequispaced data}, SIAM Journal on Scientific computing, 14 (1993), pp.~1368--1393.

\bibitem{genton2002nonparametric}
{\sc M.~G. Genton and D.~J. Gorsich}, {\em Nonparametric variogram and covariogram estimation with {F}ourier--{B}essel matrices}, Computational Statistics \& Data Analysis, 41 (2002), pp.~47--57.

\bibitem{greengard2004accelerating}
{\sc L.~Greengard and J.-Y. Lee}, {\em Accelerating the nonuniform fast {F}ourier transform}, SIAM review, 46 (2004), pp.~443--454.

\bibitem{greengard2007fast}
{\sc L.~Greengard, J.-Y. Lee, and S.~Inati}, {\em The fast sinc transform and image reconstruction from nonuniform samples in k-space}, Communications in Applied Mathematics and Computational Science, 1 (2007), pp.~121--131.

\bibitem{hansen1985fast}
{\sc E.~Hansen}, {\em Fast {H}ankel transform algorithm}, IEEE transactions on acoustics, speech, and signal processing, 33 (1985), pp.~666--671.

\bibitem{higgins1988hankel}
{\sc W.~E. Higgins and D.~C. Munson}, {\em A {H}ankel transform approach to tomographic image reconstruction}, IEEE transactions on medical imaging, 7 (1988), pp.~59--72.

\bibitem{jiang2023dual}
{\sc S.~Jiang and L.~Greengard}, {\em A dual-space multilevel kernel-splitting framework for discrete and continuous convolution}, arXiv preprint arXiv:2308.00292,  (2023).

\bibitem{johansen1979fast}
{\sc H.~Johansen and K.~S{\o}rensen}, {\em Fast {H}ankel transforms}, Geophysical Prospecting, 27 (1979), pp.~876--901.

\bibitem{johnson1987improved}
{\sc H.~F. Johnson}, {\em An improved method for computing a discrete hankel transform}, Computer physics communications, 43 (1987), pp.~181--202.

\bibitem{kapur1995algorithm}
{\sc S.~Kapur and V.~Rokhlin}, {\em An algorithm for the fast {H}ankel transform}, tech. report, Technical Report 1045, Computer Science Department, Yale University, 1995.

\bibitem{li2015butterfly}
{\sc Y.~Li, H.~Yang, E.~R. Martin, K.~L. Ho, and L.~Ying}, {\em Butterfly factorization}, Multiscale Modeling \& Simulation, 13 (2015), pp.~714--732.

\bibitem{linton2006schlomilch}
{\sc C.~Linton}, {\em Schl{\"o}milch series that arise in diffraction theory and their efficient computation}, Journal of Physics A: Mathematical and General, 39 (2006), p.~3325.

\bibitem{liu1999nonuniform}
{\sc Q.~H. Liu and Z.~Q. Zhang}, {\em Nonuniform fast hankel transform ({NUFHT}) algorithm}, Applied optics, 38 (1999), pp.~6705--6708.

\bibitem{lord1954a}
{\sc R.~D. Lord}, {\em The use of the {H}ankel transform in statistics {I}. {G}eneral theory and examples}, Biometrika, 41 (1954), pp.~44--55.

\bibitem{lord1954b}
{\sc R.~D. Lord}, {\em The use of the {H}ankel transform in statistics {II}. {M}ethods of computation}, Biometrika, 41 (1954), pp.~344--350.

\bibitem{marshall2023fast}
{\sc N.~F. Marshall, O.~Mickelin, and A.~Singer}, {\em Fast expansion into harmonics on the disk: {A} steerable basis with fast radial convolutions}, SIAM Journal on Scientific Computing, 45 (2023), pp.~A2431--A2457.

\bibitem{mook1983algorithm}
{\sc D.~Mook}, {\em An algorithm for the numerical evaluation of the {H}ankel and {A}bel transforms}, IEEE transactions on acoustics, speech, and signal processing, 31 (1983), pp.~979--985.

\bibitem{nochetto2009theory}
{\sc R.~H. Nochetto, K.~G. Siebert, and A.~Veeser}, {\em Theory of adaptive finite element methods: an introduction}, in Multiscale, Nonlinear and Adaptive Approximation: Dedicated to Wolfgang Dahmen on the Occasion of his 60th Birthday, Springer, 2009, pp.~409--542.

\bibitem{nolan2013multivariate}
{\sc J.~P. Nolan}, {\em Multivariate elliptically contoured stable distributions: theory and estimation}, Computational statistics, 28 (2013), pp.~2067--2089.

\bibitem{olver2010nist}
{\sc F.~W. Olver}, {\em NIST handbook of mathematical functions}, Cambridge university press, 2010.

\bibitem{oneil2010algorithm}
{\sc M.~O'Neil, F.~Woolfe, and V.~Rokhlin}, {\em An algorithm for the rapid evaluation of special function transforms}, Applied and Computational Harmonic Analysis, 28 (2010), pp.~203--226.

\bibitem{oppenheim1980computation}
{\sc A.~V. Oppenheim, G.~V. Frisk, and D.~R. Martinez}, {\em Computation of the {H}ankel transform using projections}, The Journal of the Acoustical Society of America, 68 (1980), pp.~523--529.

\bibitem{pang2020interpolative}
{\sc Q.~Pang, K.~L. Ho, and H.~Yang}, {\em Interpolative decomposition butterfly factorization}, SIAM Journal on Scientific Computing, 42 (2020), pp.~A1097--A1115.

\bibitem{rangan2020factorization}
{\sc A.~Rangan, M.~Spivak, J.~And{\'e}n, and A.~Barnett}, {\em Factorization of the translation kernel for fast rigid image alignment}, Inverse Problems, 36 (2020), p.~024001.

\bibitem{siegman1977quasi}
{\sc A.~Siegman}, {\em Quasi fast {H}ankel transform}, Optics letters, 1 (1977), pp.~13--15.

\bibitem{thakur2011synchrosqueezing}
{\sc G.~Thakur and H.-T. Wu}, {\em Synchrosqueezing-based recovery of instantaneous frequency from nonuniform samples}, SIAM Journal on Mathematical Analysis, 43 (2011), pp.~2078--2095.

\bibitem{townsend2015fast}
{\sc A.~Townsend}, {\em A fast analysis-based discrete {H}ankel transform using asymptotic expansions}, SIAM Journal on Numerical Analysis, 53 (2015), pp.~1897--1917.

\bibitem{watson1922treatise}
{\sc G.~N. Watson}, {\em A treatise on the theory of Bessel functions}, vol.~2, The University Press, 1922.

\bibitem{wimp1962polynomial}
{\sc J.~Wimp}, {\em Polynomial expansions of {B}essel functions and some associated functions}, Mathematics of Computation, 16 (1962), pp.~446--458.

\bibitem{zhao2013fourier}
{\sc Z.~Zhao and A.~Singer}, {\em Fourier--{B}essel rotational invariant eigenimages}, JOSA A, 30 (2013), pp.~871--877.

\bibitem{zhou2022spectral}
{\sc R.~Zhou and N.~Grisouard}, {\em Spectral solver for {C}auchy problems in polar coordinates using discrete {H}ankel transforms}, arXiv preprint arXiv:2210.09736,  (2022).

\end{thebibliography}

\end{document}